\documentclass[11pt,twoside, leqno]{article}

\usepackage{amssymb,amsmath,amsfonts,amsthm,color,mathrsfs}
\usepackage[Symbol]{upgreek}
\usepackage{txfonts}
\usepackage[nottoc,notlot,notlof]{tocbibind}
\usepackage[active]{srcltx}
\usepackage[dvipsnames, svgnames, x11names]{xcolor}

\usepackage{graphicx}
\usepackage{epsfig}

\usepackage{txfonts}
\allowdisplaybreaks
\def\ncc{{\mathcal C}}
\def\ncw{{\mathcal W}}

\def\L{{\mathcal L}}
\def\K{\mathbf{K}}
\def\rr{{\mathbb R}}
\def\rn{{{\rr}^n}}

\def\cn{{\mathbb N}}

\def\cs{{\mathcal S}}

\def\cf{{\mathcal F}}
\def\cg{{\mathcal G}}

\def\ca{{\mathcal A}}

\def\supp{{\mathop\mathrm{\,supp\,}}}
\def\dist{{\mathop\mathrm {\,dist\,}}}
\def\diam{{\mathop\mathrm {\,diam\,}}}

\def\r{\right}
\def\lf{\left}

\newtheorem{thm}{Theorem}[section]
\newtheorem{lem}[thm]{Lemma}
\newtheorem{prop}[thm]{Proposition}
\newtheorem{cor}[thm]{Corollary}

\newtheorem{conje}[thm]{Conjecture}

\newtheorem{rem}[thm]{Remark}

\numberwithin{equation}{section}

\begin{document}
\arraycolsep=1pt
\author{Ren-Jin Jiang,  Hong-Quan Li, Hai-Bo Lin}
\title{{\bf Riesz transform on manifolds with ends of different volume growth for $1<p<2$}
 \footnotetext{\hspace{-0.35cm} 2010 {\it Mathematics
Subject Classification}. Primary  42B20; Secondary  58J35, 43A85.
\endgraf{
{\it Key words and phrases: Riesz transform, heat kernel, bounded geometry, non-doubling measure}
\endgraf}}
\date{}}
\maketitle

\begin{center}
\begin{minipage}{11.5cm}\small
{\noindent{\bf Abstract}.
Let $M_1$, $\cdots$, $M_\ell$ be complete, connected and non-collapsed  manifolds of the same dimension, where $2\le \ell\in\mathbb{N}$, and suppose that each $M_i$ satisfies  a doubling condition and a Gaussian upper bound for the heat kernel.
If each manifold $M_i$ has volume growth either bigger than two or equal to two, then we show that the Riesz transform $\nabla \L^{-1/2}$
is bounded on $L^p(M)$ for each $1<p<2$ on the gluing manifold $M=M_1\#M_2\#\cdots \# M_\ell$.
}\end{minipage}
\end{center}
\vspace{0.2cm}
\tableofcontents

\section{Introduction}

\hskip\parindent In this paper, we consider a complete, connected and  non-compact Riemannian manifold
$M$, that is obtained by gluing  together several complete manifolds of the same dimension.
On the manifold $M$, we denote by $d$ the geodesic distance, by $\mu$ the
Riemannian measure. We denote by $\L$  the non-negative Laplace-Beltrami operator   on $M$, and  let $\{e^{-t\mathscr{L}}\}_{t>0}$ be the heat semigroup.
 The corresponding Riesz transform is then given by
  $$\nabla\L^{-1/2}=\frac{1}{\sqrt{\pi}} \int_0^\infty \nabla e^{-s\L}\frac{\,ds}{\sqrt s},$$
where $\nabla$ denotes the Riemannian gradient. For study and developments on the Riesz transform, we refer the readers to \cite{Al92,An92,ac05,acdh,bak2,bf15,bk04,ca07,ca16,cch06,ch92,CMO15,cd99,cd03,cjks16,
lh99,ls21,lx10,lohoue85,lohoue92,shz05,str83,var88} and references therein.
Notice that as a consequence of integration by parts, the Riesz transform $\nabla\L^{-1/2}$ is always bounded on $L^2(M)$. { We emphasize that the classical result on Euclidean spaces is not always valid in this general setting. More precisely, for any fixed $q_0 \ge 2$, there exists Riemannian manifolds of this type where the Riesz transform is bounded on $L^p$ for $1 < p < q_0$ but unbounded for $p > q_0$; if $q_0 > 2$ it is unbounded also on $L^{q_0}$ and not even of weak type $(q_0, \, q_0)$. Cf. e.g.\cite{cd99, lh99} for early concrete examples.}

To move further, let us recall some basic notation. We denote by $B(x,r)$, $B_i(x_i,r)$ the open ball with centre $x\in M$, $x_i\in M_i$ and radius $r>0$ in $M$, $M_i$, and by $V(x,r)$, $V_i(x_i,r)$ their volume $\mu(B(x,r))$, $\mu_i(B_i(x_i,r))$, respectively.
We say that $M_i$ satisfies the volume doubling property (in short is  doubling)  if there exists  a constant  $C_{D}>1$ such that
$$ V_i(x_i,2r)\le C_{D}V_i(x_i,r), \qquad \forall r > 0, \  x_i\in M_i. \leqno(D)$$
Let $\mathcal{L}_i$ denote the non-negative Laplace-Beltrami operator on $M_i$.
The heat semigroup on $M_i$ has a smooth positive and symmetric kernel $h_{i,t}(x,y)$,  meaning that
$$e^{-t\mathcal{L}_i}f(x)=\int_{M_i} h_{i,t}(x,y)f(y)\,d\mu_i(y)$$
for suitable  functions $f$.
One  says that the heat kernel satisfies a Gaussian upper  bound if
there exists a constant $C >0$ such that
$$h_{i,t}(x,y)\le
  \frac{C}{V_i(x,{\sqrt t})}\exp\lf\{-\frac{d^2(x,y)}{C t}\r\}, \qquad \forall t > 0, \ x, y\in M_i.\leqno(UE)
$$

By Coulhon and Duong \cite{cd99}, the Riesz transform is bounded on $L^p(M_i)$ for all $p\in (1,2)$ if
$(D)$ and $(UE)$ are satisfied. Chen et al. \cite{CCFR} and Li-Zhu \cite{lz17} showed further that
the Gaussian upper bound can be relaxed.
We note that for a closely related operator, the Littlewood-Paley-Stein operator, Coulhon, Duong and Li
in \cite{cdl03} proved that the  operator is always $L^p$-bounded for $1<p<2$ on complete manifolds. Meanwhile, Coulhon-Duong \cite{cd03} raised the following conjecture.
\begin{conje}[Coulhon-Duong]
The Riesz transform is $L^p$-bounded for $1<p<2$ on complete manifolds.
\end{conje}
The conjecture is far from being proved or disproved  to the best of our knowledge, as there are no efficient tools to deal with the Riesz transform in general setting.

In this paper, we shall focus on $L^p$-boundedness of $\nabla\L^{-1/2}$ for $1<p<2$.
Let us consider a gluing manifold $M$, given as
$$M=M_1\#M_2\#\cdots \# M_\ell=E_0 \cup \cup_{i=1}^\ell E_i,$$
where each $M_i$ is a complete manifold and all $M_i$ are of the same dimension, $E_i=M_i\setminus K_i$, $K_i$ and $E_0$ are compact sub-manifolds; see Figure 1. We refer the reader to \cite[\S~2.2]{gri-sal09} for further details.

\begin{figure}[ht]
\centerline{ \epsfig{file=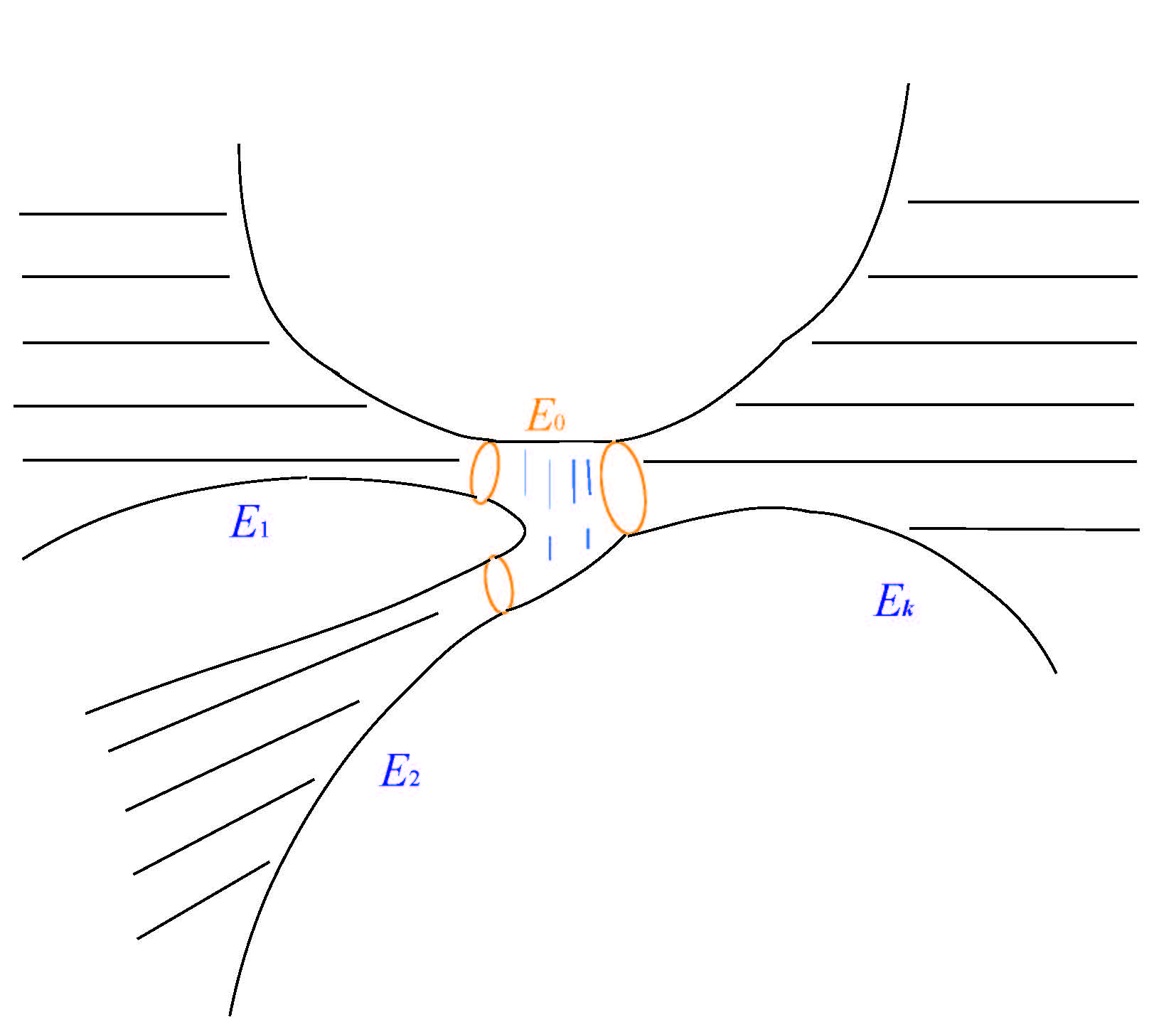, scale=0.5}
              }
             \caption{A gluing manifold.}
\end{figure}

If the gluing manifold $M$ also satisfies $(D)$, then $(UE)$ is preserved by the gluing result of Grigor'yan and Saloff-Coste \cite{gri-sal09}, and the Riesz transform is bounded on $L^p(M)$ for $1<p<2$ by Coulhon and Duong \cite{cd99}; see also Sikora \cite{Si}.

We are interested in cases where the doubling condition fails on $M$. Let us review some closely related results.
Hassell and Sikora \cite[Theorem 7.1]{hs19} proved the following. Suppose that for each $i$,
$M_i$ has lower Ricci curvature bound and positive injectivity radius. Moreover, suppose $M_i$ satisfies $(D)$ and it holds that
\begin{eqnarray}\label{hs-1}
V_i(x,r)\le \begin{cases}
Cr^{N} & r\le 1, \\
Cr^{n_i} & r>1,
\end{cases}, \qquad \forall x \in M_i,
\end{eqnarray}
and
\begin{eqnarray}\label{hs-2}
\left\|\sqrt t \nabla e^{-t\mathcal{L}_i}\right\|_{1\to\infty}
\le \begin{cases}
Ct^{-N/2} & t\le 1, \\
Ct^{-n_i/2} & t>1,
\end{cases}
\end{eqnarray}
where $n_i\ge 3$, $\|\cdot\|_{1\to\infty}$ denotes the operator norm from $L^1$ to $L^\infty$. Then the Riesz transform is weakly $(1,1)$ bounded, and $(p,p)$ bounded iff $1<p<\min_i\{n_i\}$.
Very recently, Hassell, Nix and Sikora \cite{hns19} further obtained weak $(1,1)$ boundness
on the same type manifold but allowing $n_i\ge 2$.

\begin{rem}\rm
As pointed out in \cite{hs19} (see (59) and the proof of Theorem 7.1 there), it follows from  \cite[Proposition 2.1 \& Corollary 2.2]{CS10} that,
the conditions \eqref{hs-1} and \eqref{hs-2} together imply that,
the heat kernel $h_{i,t}(x,y)$ of  $e^{-t\mathcal{L}_i}$ satisfies the Li-Yau estimate, i.e.,
$$ \frac{1}{C V_i(x,{\sqrt t})}\exp\lf\{- C \frac{d^2(x,y)}{t}\r\}\le h_{i,t}(x,y)\le
  \frac{C}{V_i(x,{\sqrt t})}\exp\lf\{-\frac{d^2(x,y)}{C t}\r\},  \leqno(LY)
$$
and the volume of balls has growth as
\begin{eqnarray}\label{volume-ahlfors}
C^{-1} r^N \le V_i(x,r) \le C r^N, \  \forall \, 0 < r \le 1, \quad C^{-1} r^{n_i} \le V_i(x,r) \le C r^{n_i}, \  \forall \, r \ge 1,
\end{eqnarray}
\end{rem}

Previous to \cite{hs19,hns19}, Carron \cite{ca07} obtained some stability result
for $p\in (\frac{\nu}{\nu-1},\nu)$, where $\nu>3$ is from the global Sobolev inequality.
Devyver \cite{de15} further refined the Sobolev dimension $\nu$ to hyperbolic dimension.

It is then nature to wonder whether one can prove $L^p$-boundedness of the Riesz transform for all $1<p<2$ without requiring polynomial growth of the volume \eqref{volume-ahlfors}, or the gradient estimates of heat kernel \eqref{hs-2}?
Recall that a manifold $\mathrm{M}$ is said to be non-collapsed, provided that
the volume each ball with radius one in $\mathrm{M}$ has a positive bottom, i.e., $\inf_{x\in \mathrm{M}}V(x,1)>0$. In the paper, we shall prove the following:
\begin{thm}\label{main-result}
Let $2 \le \ell \in \cn$. Suppose that $\{M_i\}_{i=1}^\ell$ are complete, connected and  non-collapsed manifolds of the same dimension, and each $M_i$ satisfies  $(D)$ and $(UE)$. Moreover, assume that there exist constants $C, c > 0$ and some $x_i\in M_i$ for each $i$ such that for all $R\ge r\ge 1$, it holds either
$$c \left(\frac{R}{r}\right)^{2} \le \frac{V_{i}(x_i,R)}{V_{i}(x_i,r)} \le C \left(\frac{R}{r}\right)^{2},$$
or for some $n_i>2$
$$c \left(\frac{R}{r}\right)^{n_i} \le \frac{V_{i}(x_i,R)}{V_{i}(x_i,r)}.$$
Let $M=M_1\#M_2\#\cdots \# M_\ell$. Then the Riesz transform $\nabla \L^{-1/2}$
is bounded on $L^p(M)$ for each $1<p<2$.
\end{thm}

\begin{rem}\rm
(i) We do not know how to prove weak $(1,1)$ boundedness. Our method uses the mapping property of the operators $\nabla e^{-t\L}$ (see Proposition \ref{mapping-gradient-heat} below)
and $t\L e^{-t\L}$ (see Proposition \ref{time-heat-map}), which only have optimal bounds for $p>1$.

(ii) The requirement of the case $n_i=2$ is stronger than the case $n_i>2$, since $n_i=2$ corresponds to the critical case, the behavior of the heat kernel is much more complicate and will be not
handful if we do not assume the upper bound of the volume growth; see \cite{gri-sal09}. The case $n_i<2$
is missing from the main result due to the same reason. In fact, by assuming (stronger) two-side Gaussian bounds of the heat kernel on each $M_i$, \cite[Section 6]{gri-sal09} gives optimal heat kernel estimates including the cases that some $n_i<2$, which we believe can be used to prove $L^p$-boundedness of the Riesz transform for $1<p<2$. However, since the present proof is already quite long, we will deal these cases in future.

(iii) The case that all $n_i$ equal two is obvious, since in this case, the manifold is doubling and the
heat kernel satisfies $(UE)$ by \cite{gri-sal09}, the required conclusion follows from
\cite{cd99}. Therefore, we will only prove the case that some $n_i>2$.

(iv) Since we only need a upper Gaussian bound of the heat kernel on $M_i$, our result applies to any
uniformly elliptic operators on these manifolds (\eqref{hs-2} in general does not hold). Moreover, our result can be applied to the case
where the manifolds have volume growth different from \eqref{volume-ahlfors}, which seems to be also new (see \cite[Remark 7.3]{hs19}).
Let us take an example from \cite{gis18}. For $\alpha\in (0,2)$, consider $\mathcal{R}^\alpha:=(\rr^2,g_\alpha)$, where $g_\alpha$
is a Riemannian metric such that, in the polar coordinates $(\rho,\theta)$, for $\rho>1$ it equals
$$g_\alpha=\,d\rho^2+\rho^{2(\alpha-1)}\,d\theta^2.$$
The volume of balls $B(x,r)$ on $\mathcal{R}^\alpha$, $r>1$, has growth as
$$V(x,r)\sim\begin{cases}
r^\alpha,& \,|x|<r\\
\min\{r^2,\,r|x|^{\alpha-1}\},&\,|x|\ge r.
\end{cases}$$
In particular, $V(0,r)\sim r^\alpha$ for $r>1$. Note that for $\alpha\in (0,2)$, the exterior part $\{\rho>1\}$ of $\mathcal{R}^\alpha$
is isometric to a certain surface of revolution in $\rr^3$ and the Li-Yau estimate holds on $\mathcal{R}^\alpha$; see \cite[pp. 160-161]{gis18}.
For $n\ge 4$, let $\mathcal{M}_1$ and $\mathcal{M}_2$ be closed manifolds of dimension $n-4$
and $n-2$ respectively. Our result then applies to
the gluing manifolds $\rn\# (\rr^2\times \mathcal{R}^\alpha\times \mathcal{M}_1)$ and
$(\rr^2\times \mathcal{R}^\alpha\times \mathcal{M}_1)\# (\rr^{2}\times \mathcal{M}_2)$.
\end{rem}

We remark that the $L^p$-boundedness of $\nabla\L^{-1/2}$ for $1<p\le 2$ is best possible in our setting,
since we allow that some ends have volume growth as the plane at infinity. By
Hassell et al. \cite{hns19}, on a manifold with at least two ends, if one end has volume growth as the plane at infinity, then the Riesz transform is not bounded for any $p>2$; see also \cite{cd99} for the case $M=\rr^2\# \rr^2$ and \cite{cch06}.
In fact, even in the case that all ends are non-parabolic, i.e., $\min\{n_i\}>2$, it is known that some additional condition are needed for the case $p>2$ (cf. \cite{acdh,ji20}).

The proof of the main result is rather long and delicate, we shall explain the method in the beginning of
Section 4. Briefly speaking, we shall decompose the manifold into several subsets, and deal with each subset with different methods. The heat kernel estimate from Grigor'yan and Saloff-Coste  \cite{gri-sal09}
and the original method of Coulhon and Duong \cite{cd99}
play important roles in our proof. We shall frequently use the $L^p$-Davies-Gaffney estimates
for the operators $e^{-t\L}$, $\nabla e^{-t\L}$ and $\L e^{-t\L}$, and $L^p(E)\to L^2(F)$ type bounds
for the operators $e^{-t\L}$  and $\L e^{-t\L}$.
Note that it is {\em not} possible
to have a Gaussian  upper bound for the time derivative of the heat kernel in our setting as in \cite{cd99}; see \cite{gri-sal09}.
We shall use a powerful tool established by Davies \cite{davies97} to get an upper bound for
the time derivative of the heat kernel $\partial_th_t(x,y)$, which provides us the right norm bound
of the operator $\L e^{-t\L}$ from $L^p(E)$ to $L^2(F)$ for some subsets $E,\,F$. We believe this is one of  the main ingredients of the paper, see Section 3 and Subsection 4.4.

The paper is organized as follows, In Section 2, we provide some basic estimates for the heat kernel, and the local Riesz transform. In Section 3, we provide estimates and mapping property for the heat kernel and its time derivative. In Section 4, we provide the proof for the main result.

Throughout the work, we denote by $C,c$ positive constants which are independent of the
main parameters, but which may vary from line to line. The symbol $A \lesssim B$ means
that $A\le CB$, $A\lesssim_{\alpha,\beta}B$ means that the implicit constant depends on $\alpha,\beta$,  and $A\sim B$ means $cA\le B\le CA$, for some harmless constants $c,C>0$.
We use $\|\cdot\|_{p}$ to denote the $L^p(M)$ norm, and
 $\|\cdot\|_{p\to q}$ to denote the operator norm $\|\cdot\|_{L^p(M)\to L^q(M)}$.

\section{Preliminaries}
\hskip\parindent In this section, we collect some basic estimates of the heat semigroup
and heat kernels.
\subsection{Davies-Gaffney estimates}
\hskip\parindent In this part, we provide some basic estimate for the heat semigroup. The first result is in fact a direct consequence of \cite[Theorem 4.1]{cd03} and \cite[Proposition 3.6]{CS10}:
\begin{prop}\label{mapping-gradient-heat}
The operator $\nabla e^{-t\L}$ is bounded on $L^p(M)$ for $1<p\le 2$ with
$$\|\nabla e^{-t\L}\|_{p\to p}\lesssim_p \frac{1}{\sqrt t}, \qquad \forall t > 0.$$
\end{prop}
\begin{proof}By \cite[Theorem 4.1]{cd03}, for any $f\in C^\infty_0(M)$ and $p\in (1,2]$, we have
$$\|\nabla e^{-t\L}f\|_{p}^2\lesssim_p \|e^{-t\L}f\|_{p} \, \|\L e^{-t\L}f\|_{p}\lesssim_p \frac{1}{t}\|f\|_{p}^2,$$
by the classical Littlewood-Paley-Stein theory (c.f. \cite{S70}).
\end{proof}

For $1\le p<\infty$, we say that an operator $T$ satisfies the $L^p$-Davies-Gaffney estimate, if there exists $C>0$ such that
for any closed sets $E,F\subset M$, it holds
\begin{eqnarray*}
\|T(f\chi_E)\|_{L^p(F)}\le C\exp\left(-\frac{d(E,F)^2}{C t}\right)\|f\|_{L^p(E)}.
\end{eqnarray*}
When $p=2$, we shall say that $T$ satisfies the Davies-Gaffney estimate for short.

The following result was proved in \cite[p. 930, (3.1)]{acdh}; see also \cite{CS}.

\begin{prop}\label{davies-gaffney-off}
The operators  $e^{-t\L}$, $\sqrt t\nabla e^{-t\L}$ and $t\L e^{-t\L}$ satisfy the Davies-Gaffney estimate.
\end{prop}

Using the Riesz-Thorin interpolation theorem, we deduce the following $L^p$-Davies-Gaffney estimates.
\begin{cor}\label{davies-operators}
The operators $e^{-t\L}$, $\sqrt t\nabla e^{-t\L}$ and $t\L e^{-t\L}$  satisfy $L^p$-Davies-Gaffney estimate for $1<p\le 2$.
\end{cor}
\begin{proof}
Note that $e^{-t\L}$  and $t\L e^{-t\L}$ are $L^p$-bounded for all $1<p<\infty$,
and by Proposition \ref{mapping-gradient-heat}, $\sqrt{t}\nabla e^{-t\L}$ is bounded on $L^p$ for $1<p\le 2$.
From this, together with the Riesz-Thorin theorem and Proposition \ref{davies-gaffney-off}, we deduce that for $1<p<2$,
\begin{eqnarray*}
&&\|e^{-t\L}(f\chi_E)\|_{L^p(F)}+\|\sqrt t\nabla e^{-t\L}(f\chi_E)\|_{L^p(F)}+\|t\L e^{-t\L}(f\chi_E)\|_{L^p(F)}\\
&&\quad \le
C(p)\exp\left(-c(p)\frac{d(E,F)^2}{t}\right)\|f\|_{L^p(E)},
\end{eqnarray*}
as desired.
\end{proof}
From Corollary \ref{davies-operators} and \cite[Proposition 3.1]{au07}, it follows that the
composition of $\sqrt t\nabla e^{-t\L}$ and $t\L e^{-t\L}$ satisfies:
\begin{cor}\label{davies-gaffney-com}
Let $1<p\le 2$. The operator $t^{3/2}\nabla \L e^{-t\L}$ satisfies the $L^p$-Davies-Gaffney estimate.
\end{cor}

\subsection{Notation and basic properties on volume growth}
\hskip\parindent
We will need the heat kernel estimates obtained in  \cite{gri-sal09}. For this purpose, let us begin with some basic properties of the volume growth.

As in \cite[\S~4.3]{gri-sal09}, let $B(x, r)$ (resp. $B_i(x, r)$ with $1 \leq i \leq \ell$) denote the geodesic ball in $M$ (resp. $M_i$). We set
\begin{align}
|x|:=\sup_{y\in E_0}\{d(x,y)\}, \quad V(x, r) :=V(B(x, r)),  \quad   V_i(x,r):=V_i(B_i(x,r)),  \quad  V_i(r) :=V_i(o_i,r),
\end{align}
 where  $o_i\in \partial E_i$ is a fixed reference point. Note that  $\partial E_i$ is the set that connecting $E_i=M_i\setminus K_i$ to $E_0$, and it always holds that
\begin{equation}\label{com-metric}
d(x,y)\le |x|+|y|.
\end{equation}
In what follows, for simplicity of notions, we shall assume that
\begin{equation}\label{c1}
\mbox{diam}(E_0)=1,\, \mu(E_0)=1, \,  1\le V_i(1)=V_i(o_i,1) \le 4, \, \forall\, 1 \le i \le \ell,
\end{equation}
and
\begin{align}
1 \leq \mu(F_*) \leq 40 \ell, \  \mbox{with} \  F_* := \{x \in M; \, \dist(x, E_0) \leq 2\}.
\end{align}

Note that the doubling condition on $M_i$ ($1 \leq i \leq \ell$) implies that,
 there exists $N_i$, which is of course not less than $n_i$, such that
\begin{align} \label{c2}
\frac{V_i(x,R)}{V_i(x,r)}\lesssim  \left(\frac{R}{r}\right)^{N_i}, \qquad \forall x \in M_i, \  0 < r \le R.
\end{align}
Throughout the paper, we set
\begin{equation}\label{set-index}
 N_\infty:=\max_i\left\{N_i:\,i=1,\cdots,\ell\right\},
\end{equation}

We also set
\begin{align} \label{na}
V_0(r):=\inf_{1\le i\le \ell} V_i(r), \ r > 0; \qquad F_i^{(r)} := \{x \in E_i; \, \dist(x, E_0) \leq 2 r\}, \ r \ge 1, \ 1 \le i \le \ell.
\end{align}
And the following volume growth properties will be used:

\begin{lem} \label{vg}
(i) It holds that
\begin{align} \label{vg1}
\frac{V_i(R)}{V_i(r)} \gtrsim \left( \frac{R}{r} \right)^{n_i} \ge \left( \frac{R}{r} \right)^2, \qquad \forall 1 \le r \le R, \  1 \le i \le \ell.
\end{align}
In particular, we have $V_0(r) \gtrsim r^2$ for all $r \ge 1$.

(ii) We have
\begin{align} \label{vcp}
\mu(F_i^{(r)}) \sim V_i(2 r) \sim V_i(r), \qquad \forall r \ge 1, \  1 \le i \le \ell.
\end{align}

(iii) $M$ satisfies the local doubling volume property
\begin{align*}
\frac{V(x, 2 r)}{V(x, r)} \lesssim 1, \qquad \forall x \in M, \  0 < r \le 1,
\end{align*}
and the volume is of at most polynomial growth in the sense that
\begin{align*}
\frac{V(x, r)}{V(x, 1)} \lesssim r^{N_{\infty}}, \qquad \forall x \in M, \  r \ge 1.
\end{align*}
\end{lem}

\begin{proof}
We start by proving (i). Indeed for each $M_i$, recall that there exists $x_i \in M_i$ such that
\begin{align*}
\frac{V_i(x_i,R)}{V_i(x_i,r)} \gtrsim \left( \frac{R}{r} \right)^{n_i}, \qquad \forall 1 \le r \le R.
\end{align*}
By the simple fact that
\begin{align*}
B_i(o_i, R) \subset B_i\big(x_i, 2 (R+d(x_i, o_i)) \big) \subset B_i\big(o_i, 4 (R+d(x_i, o_i)) \big), \quad B_i(o_i, r) \subset B_i\big(x_i, 2 (r+d(x_i, o_i)) \big),
\end{align*}
we can write
\begin{align*}
\frac{V_i(R)}{V_i(r)} &\ge \frac{V_i(R)}{V_i(x_i, 2 (R + d(x_i, o_i)))} \, \frac{V_i(x_i, 2 (R + d(x_i, o_i)))}{V_i(x_i, 2 (r + d(x_i, o_i)))}  \\
&\ge \frac{V_i(R)}{V_i(4 (R + d(x_i, o_i)))} \, \frac{V_i(x_i, 2 (R + d(x_i, o_i)))}{V_i(x_i, 2 (r + d(x_i, o_i)))}.
\end{align*}
Then \eqref{c2} implies that
\begin{eqnarray*}
\frac{V_i(R)}{V_i(r)}
&&\gtrsim \left(\frac{R}{R+d(x_i,o_i)}\right)^{N_i} \, \left(\frac{R+d(x_i,o_i)}{r+d(x_i,o_i)}\right)^{n_i}
= \left(\frac{R}{r}\right)^{n_i}\left(\frac{R}{R+d(x_i,o_i)}\right)^{N_i-n_i}
\left(\frac{r}{r+d(x_i,o_i)}\right)^{n_i}\\
&&\gtrsim \left(\frac{R}{r}\right)^{n_i}\left(\frac{1}{1+d(x_i,o_i)}\right)^{N_i-n_i}
\left(\frac{1}{1+d(x_i,o_i)}\right)^{n_i}\\
&&\gtrsim \left(\frac{R}{r}\right)^{n_i},
\end{eqnarray*}
which completes the proof.

Using \eqref{c2} and the simple fact that $B_i(o_i, 2 r) \subset F_i^{(r)}\cup K_i \subset B_i(o_i, 4 r)$ for all $r \ge 1$, the claims (ii) and (iii) are clear.
\end{proof}

\subsection{Heat kernel upper bounds}
\hskip\parindent
Let
$$H(x,t):=\min\left\{1; \frac{|x|^2}{V_{i_x}(|x|)}+\left(\int_{|x|^2}^t\frac{\,d s}{V_{i_x}(\sqrt s)}\right)_+\right\},$$
where $(\cdot)_+$ denotes the non-negative part and $i_x$ denotes the index of the end that $x$ belongs to.

Under our assumptions of the volume growth and the doubling condition,
$$H(x,t) \sim \begin{cases}
\frac{|x|^2}{V_{i_x}(|x|)},& \, n_{i_x}>2\\
1,&\, n_{i_x}=2
\end{cases}$$
(see  also  \cite[(4.21)]{gri-sal09} for the case where $n_{i_x}>2$). Therefore, we have the uniform bound
\begin{align} \label{h}
H(x,t) \lesssim \frac{|x|^2}{V_{i_x}(|x|)} \lesssim 1.
\end{align}

We will use repeatedly the following heat kernel upper bounds on $M$, which can be deduced from \cite[\S~4]{gri-sal09}, especially \cite[Theorem 4.9]{gri-sal09} and its remarks therein:

\begin{thm}[Grigor'yan \& Saloff-Coste] \label{hke}
(i) The small time heat kernel Gaussian upper bounds hold, namely
\begin{align} \label{stu}
h_t(x,y) \lesssim  \frac{1}{V(x,\sqrt t)}\exp\left(- c \frac{d(x,y)^2}{t}\right), \qquad \forall\, 0< t \le 1, \ x, y \in M.
\end{align}
(ii)
We have that
\begin{align} \label{HEa}
h_t(x, y) &\lesssim
\frac{1}{V_0(\sqrt{t})} \frac{|x|^2}{V_i(|x|)}  \frac{|y|^2}{V_i(|y|)}
\, e^{-c \frac{|x|^2 + |y|^2}{t}} \nonumber \\
&\quad+ \min\bigg(\frac{1}{V_i(x, \sqrt{t})},  \frac{1}{V_i(y, \sqrt{t})} \bigg) \, \exp\bigg\{ -c \frac{d(x, y)^2}{t} \bigg\},
\end{align}
for any $t\ge 1$ and $x, y \in E_i$ ($1 \leq i \leq \ell$).  \\
(iii) For $r > 0$ and $g \in E_k$ ($1 \le k \le \ell$), let
\[
V_*(g, r) := \max\{V_k(g, r), \, V_k(r) \}.
\]
It holds that
\begin{align}  \label{HEb}
h_t(x, y) \lesssim \bigg( \frac{1}{V_0(\sqrt{t})} \frac{|x|^2}{V_j(|x|)}  \frac{|y|^2}{V_i(|y|)} + \frac{1}{V_*(x, \sqrt{t})}  \frac{|y|^2}{V_i(|y|)}  + \frac{1}{V_*(y, \sqrt{t})}  \frac{|x|^2}{V_j(|x|)} \bigg) \, e^{-c \frac{|x|^2 + |y|^2}{t}},
\end{align}
for all $t \ge 1$, $x \in E_j$, $y \in E_i$ ($1 \leq i \neq j \leq \ell$).  Moreover,
\begin{align} \label{hec}
h_t(x, y) \lesssim \bigg( \frac{1}{V_0(\sqrt{t})} \frac{|x|^2}{V_j(|x|)} + \frac{1}{V_*(x, \sqrt{t})} \bigg) \, e^{-c \frac{|x|^2}{t}},
\end{align}
provided $t \ge 1$, $x \in E_j$, $y \in E_0$ ($1 \leq j \leq \ell$).
\end{thm}

\begin{proof}
The estimate \eqref{stu} comes from \cite[Corollary 4.16~0]{gri-sal09}.
The claim (ii) follows from {  \cite[(4.45)]{gri-sal09} and the symmetric property of the heat kernel.} To get \eqref{HEb} and \eqref{hec}, it suffices to use \cite[Theorem~4.9, Remark~4.10, (4.24), and Remark~4.12]{gri-sal09}.
\end{proof}

In particular, we have the following:

\begin{cor} \label{uE}
 Let $1 \le i \le \ell$, $r \ge 1$ and $x$, $y \in E_i$. Let $\beta \ge 8$, then we have that
\begin{align} \label{cEn}
h_t(x, y) \lesssim_{\beta} \min\bigg(\frac{1}{V_i(x, \sqrt{t})},  \frac{1}{V_i(y, \sqrt{t})} \bigg)  \, e^{-c \, \frac{d(x, y)^2}{t}}, \quad \forall\, 0 < t \le \beta r^2, \  d(x, E_0) + d(y, E_0) \geq 2 r.
\end{align}
\end{cor}

\begin{proof}
Using the symmetric property of the heat kernel, it suffices to establish
\begin{align*}
h_t(x, y) \lesssim_{\beta} \frac{1}{V_i(x, \sqrt{t})}  \, e^{-c \, \frac{d(x, y)^2}{t}}, \quad \forall\, 0 < t \le \beta r^2, \  d(x, E_0) + d(y, E_0) \geq 2 r.
\end{align*}

Notice that (cf. also \cite[p. 1945]{gri-sal09}) $V(x, r) = V_i(x, r)$ whenever $B_i(x, r) \subset E_i$, otherwise $V(x, r) \sim V_i(x, r)$ for all $0 < r \le 1$ and $x \in E_i$. Hence according to \eqref{stu}, we may assume that $t \ge 1$.
Without loss of generality we can suppose that $d(y, E_0) \ge r$. By \eqref{HEa}, \eqref{h}, \eqref{com-metric} and the fact that $V_0(s) \gtrsim s^2$ for all $s \ge 1$, it remains to show that
\begin{align*}
\frac{1}{t} \, \frac{|y|^2}{V_i(|y|)} \, e^{-c_1 \, \frac{|x|^2+|y|^2}{t}}  \lesssim_{c_1} \frac{1}{V_i(|y|)} \, e^{- \frac{c_1}{2} \frac{|x|^2+|y|^2}{t}} \lesssim_{c_1, \beta} \frac{1}{V_i(x, \sqrt{t})},  \qquad \forall 1 \le t \le \beta r^2, \ x \in E_i,
\end{align*}
where we have used in the first inequality the fact that $s \, e^{-s} \lesssim 1$ for any $s > 0$.

Notice that $|y| \ge d(y, o_i) \ge d(y, E_0) \ge r  \gtrsim_{\beta} \sqrt{t}$. The doubling property implies that $V_i(|y|) \gtrsim_{\beta} V_i(\sqrt{t}) = V_i(o_i, \sqrt{t})$. By the fact that $|x| \ge d(o_i, x)$, it suffices to show
\[
\frac{1}{V_i(o_i, \sqrt{t})} e^{-c \, \frac{d(o_i, x)^2}{t}}  \lesssim_c \frac{1}{V_i(x, \sqrt{t})}.
\]
Indeed, its proof is based on the standard trick of doubling property, which will be used repeatedly. More precisely,
\begin{align}\label{BE1}
\frac{1}{V_i(o_i, \sqrt t)} \, e^{- c \frac{d(o_i, x)^2}{t}} &= \frac{1}{V_i(x,\sqrt t)} \, \frac{V_i(x, \sqrt t)} {V_i(o_i, \sqrt t)} \, e^{- c \frac{d(o_i, x)^2}{t}} \nonumber \\
&\le  \frac{1}{V_i(x,\sqrt t)} \, \frac{V_i(o_i, \sqrt t + d(o_i, x))} {V_i(o_i, \sqrt t)} \, e^{- c \frac{d(o_i, x)^2}{t}} \nonumber \\
&\lesssim \frac{1}{V_i(x,\sqrt t)} \, \left( \frac{\sqrt t + d(o_i, x)}{\sqrt t} \right)^{N_{\infty}} \, e^{- c \frac{d(o_i, x)^2}{t}} \lesssim_{c} \frac{1}{V_i(x,\sqrt t)}.
\end{align}

This completes the proof of this lemma.
\end{proof}

In addition,  using again $V_0(r) \gtrsim r^2$ for all $r \ge 1$, we deduce from \cite[(4.42) and (4.45)]{gri-sal09} that the following large-time on-diagonal upper bounds of the heat kernel:

\begin{lem}
It holds that:
\begin{align} \label{Le}
h_t(x, x) \lesssim \frac{1}{t} + \left\{\begin{array}{cl}
0, & \quad x \in E_0, \\[2mm]
\frac{1}{V_i(x, \sqrt{t})}, & \quad x \in E_i \ (1 \le i \le \ell),
\end{array}
\right.
\qquad \forall t \geq \frac{1}{2}.
\end{align}
\end{lem}

\subsection{The ultracontractivity of the heat semigroup $e^{-t  \L}$ ($t \ge 1$) on $M$}

\hskip\parindent
Since each $M_i$ satisfies the non-collapsing condition, one has
$$\inf_{x\in M}V(x,1)\ge c>0.$$
Moreover, since each $M_i$ is a connected doubling manifold, there is a $\delta_i>0$ such that:
$$\left(\frac{R}{r}\right)^{\delta_i}\le \frac{V_i(x,R)}{V_i(x,r)}, \qquad \forall x \in M_i, \  0 < r \le R,$$
see \cite[p. 412]{gri09} or \cite[p. 213, Remark 8.1.15]{hkst}.
Note that $\delta_i$ might be different from $n_i$ but is not larger than $n_i$.
From this, on-diagonal upper bounds \eqref{Le} and the semigroup property imply that
$$\| e^{-t \, \L} \|_{1\to \infty}\lesssim t^{-\delta}, \quad \forall t \ge \frac{1}{2}, \
\mbox{where} \  2 \delta=\min\{2, \, \min_{1 \le i\le \ell} \delta_i\}.$$

Hence, from the Riesz-Thorin interpolation theorem and the contraction property of the heat semigroup, we get the following ultracontractivity:
\begin{equation}\label{heat-operator-norm}
\|  e^{-t \, \L} \|_{p\to q}\lesssim t^{-\delta(\frac{1}{p}-\frac 1q)}, \qquad \forall t \ge \frac{1}{2},  \  1\le p \le q \le \infty.
\end{equation}

\subsection{Local estimates}
\hskip\parindent
Recall that
$$\L^{-1/2}=\frac{1}{\sqrt{\pi}} \int_0^\infty e^{-s\L}\frac{\,ds}{\sqrt s},$$
and
$$(1+\L)^{-1/2}=\frac{1}{\sqrt{\pi}} \int_0^\infty e^{-s-s\L}\frac{\,ds}{\sqrt s}.$$
We prove in the following that, by splitting the integral in the Riesz operator into small time and
  large time, the part of the small time, i.e.,
$$\int_0^1 \nabla e^{-s\L}\frac{\,ds}{\sqrt s}$$
is bounded on $L^p(M)$ for $1<p<2$ under the assumptions of Theorem \ref{main-result}.

\begin{lem}\label{local-part}
(i) The operators
$\int_0^1 \nabla e^{-s\L}\frac{\,ds}{\sqrt s}$, $\int_1^\infty \nabla e^{-s\L}\frac{\,ds}{\sqrt s}$
are bounded on $L^2(M)$.

(ii) Under the assumptions of Theorem \ref{main-result}, the operator
$$\int_0^1 \nabla e^{-s\L}\frac{\,ds}{\sqrt s}$$
is bounded on $L^p(M)$ for $1<p<2$.
\end{lem}
\begin{proof}
First, by Proposition \ref{mapping-gradient-heat}, it follows from the semigroup property that for any $1 < p \le 2$,
\begin{eqnarray*}
\left\|\int_1^\infty \nabla e^{-s-s\L}\frac{\,ds}{\sqrt s}\right\|_{p\to p}&&\le
\int_1^\infty \left\|\nabla e^{-s-s\L}\right\|_{p\to p}\frac{\,ds}{\sqrt s}\\
&&\le \int_1^\infty e^{-s}\left\|\nabla e^{-\L}\right\|_{p\to p}\left\|e^{-(s-1)\L}\right\|_{p\to p}\frac{\,ds}{\sqrt s}\\
&&{ \lesssim} \int_1^\infty e^{-s}\frac{\,ds}{\sqrt s}\\
&&{ \lesssim 1.}
\end{eqnarray*}

Next, using the small time heat kernel Gaussian upper bound (cf. \eqref{stu}) and Lemma \ref{vg}, \cite[Theorem~1.2]{cd99} implies that the local Riesz transform $\nabla(1+\L)^{-1/2}$ is bounded on $L^p(M)$ for $1<p\le 2$. This implies that
\begin{eqnarray*}
\nabla(1+\L)^{-1/2} - \frac{1}{\sqrt \pi} \int_1^\infty \nabla e^{-s-s\L}\frac{\,ds}{\sqrt s} = \frac{1}{\sqrt \pi} \int_0^1 \nabla e^{-s-s\L}\frac{\,ds}{\sqrt s}
\end{eqnarray*}
is bounded on $L^p(M)$ for $1<p\le 2$.

On the other hand, using Proposition \ref{mapping-gradient-heat} again, we get that
\begin{eqnarray*}
\left\|\int_0^1 \nabla e^{-s-s\L}\frac{\,ds}{\sqrt s}-\int_0^1 \nabla e^{-s\L}\frac{\,ds}{\sqrt s}\right\|_{p\to p}&&=\left\|\int_0^1(1-e^{-s}) \nabla e^{-s\L}\frac{\,ds}{\sqrt s}\right\|_{p\to p}\\
&&\lesssim \int_0^1(1-e^{-s}) \frac{1}{\sqrt s}\frac{\,ds}{\sqrt s}\\
&&\lesssim 1.
\end{eqnarray*}
We conclude therefore that $\int_0^1 \nabla e^{-s\L}\frac{\,ds}{\sqrt s}$ is bounded on $L^p(M)$ { for all $1 < p \le 2$.}
Moreover, since $\nabla\L^{-1/2}$ is bounded on $L^2(M)$, we also have $\int_1^\infty \nabla e^{-s\L}\frac{\,ds}{\sqrt s}$ is bounded on $L^2(M)$.
\end{proof}

\begin{rem}\rm
From the above proof, we actually see that, if the local Riesz transform $\nabla(1+\L)^{-1/2}$
is $L^p$-bounded, $1<p<2$, then the operator $\int_0^1 \nabla e^{-s\L}\frac{\,ds}{\sqrt s}$ is bounded on $L^p(M)$.
\end{rem}

\section{Mapping property of the heat semigroup and its time derivative}
\hskip\parindent
{ Recall that for $r \ge 1$,
\[
F_i^{(r)} := \{x \in E_i:\, \dist(x, E_0) \le 2 r \}, \qquad 1 \le i \le \ell.
\]
If $\gamma$  is large enough, saying $\gamma \ge 100 \ell$,  we set in the sequels,
\begin{align}
R_i := R_i(\gamma) > 1 \  \mbox{such that} \ \mu\left(F_i^{(R_i)}\right) = \gamma.
\end{align}
It holds then from \eqref{vcp} that
\begin{align} \label{brN}
\mu(F_i^{(R_i)}) \sim V_i(R_i)\sim V_j(R_j) \sim \mu(F_j^{(R_j)}), \qquad \forall \gamma \gg 1, \ 1 \le i, j \le \ell.
\end{align}
}
See Figure 2.
\begin{figure}[ht]
\centerline{ \epsfig{file=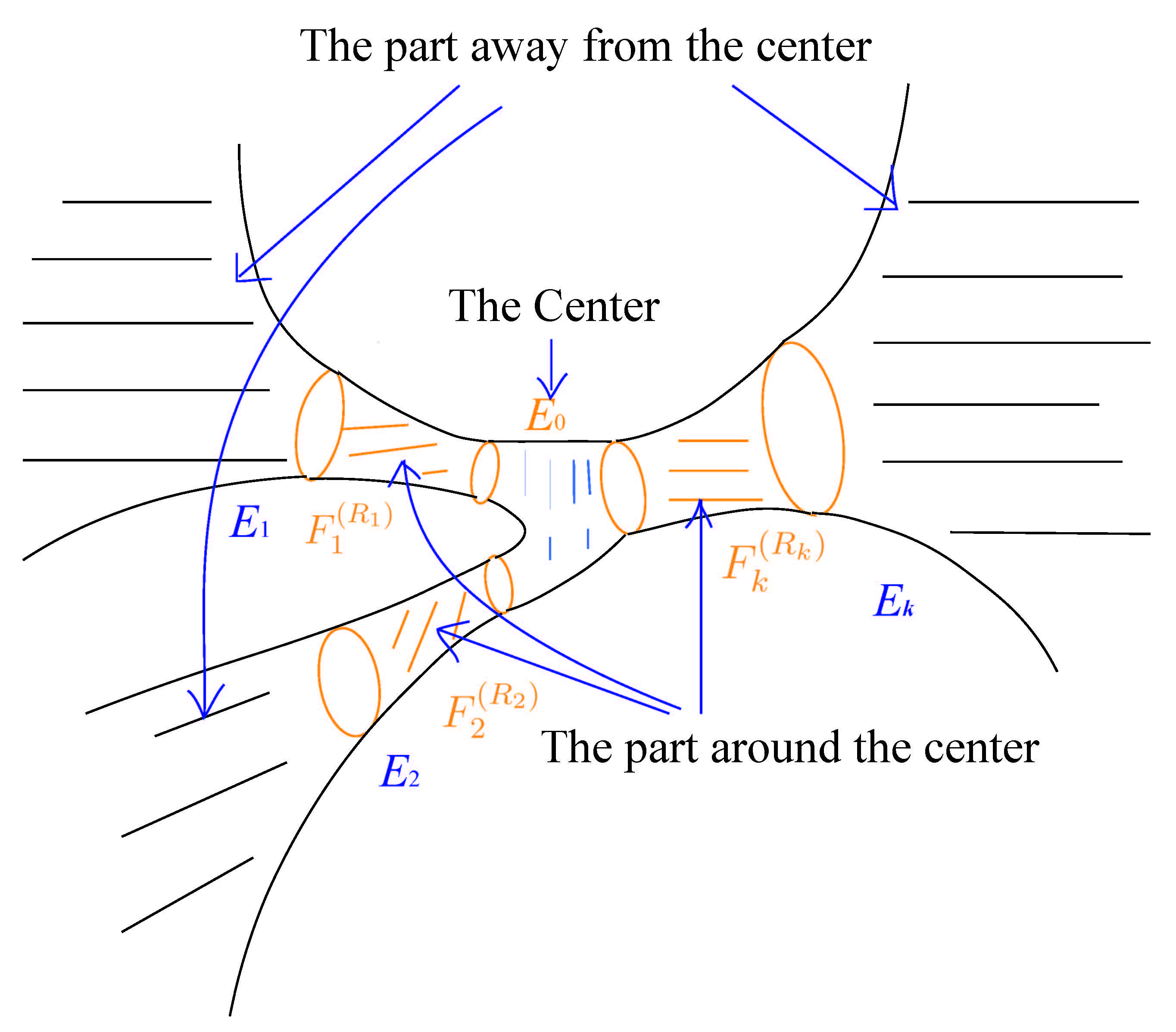, scale=0.6}
              }
             \caption{Partition of the manifold}
\end{figure}

\subsection{Mapping properties of the heat semigroup with large time}
\hskip\parindent
The following
proposition
follows directly from Theorem \ref{hke}.

\begin{prop}\label{heat-map} Under the assumptions of Theorem \ref{main-result},  for each $p\in [1,2)$, it holds
that:
$$\left\| e^{-t\L} \right\|_{L^p(F_i^{(R_i)})\to L^2(E_j\setminus F_j^{(R_j)})} \lesssim
\frac{R_j}{\sqrt{t}} \mu(F_i^{(R_i)})^{\frac 12-\frac 1p}, \qquad \forall t \ge R_j^2 = R_j(\gamma)^2, \  \gamma \gg 1,$$
and
$$\left\| e^{-t\L} \right\|_{L^p(F_i^{(R_i)})\to L^\infty(E_j\setminus F_j^{(R_j)})} \lesssim
\frac{R_j^{2}}{t}\mu(F_i^{(R_i)})^{-\frac 1p}, \qquad \forall t \ge R_j^2= R_j(\gamma)^2, \  \gamma \gg 1.$$
\end{prop}
\begin{proof}
Recall that (cf. Lemma \ref{vg} (i))
\begin{align*}
V_k(r) \ge V_0(r) \gtrsim r^2, \  \forall r \ge 1, \qquad \frac{|g|^2}{V_k(|g|)} \lesssim 1, \  \forall g \in E_k \ (1 \le k \le \ell).
\end{align*}

We start with the case $i = j$.  According to {  \eqref{HEa}, we have that
\begin{eqnarray} \label{Pe1}
h_t(x, y) &&\lesssim   \frac{1}{t} \frac{|x|^2}{V_i(|x|)}
+ \frac{1}{V_i(y, \sqrt t)} \nonumber \\
&&= \frac{1}{t} \, \frac{s^{2}}{V_i(s)} \, \bigg\{ \bigg[ \frac{V_i(|x|)}{V_i(s)} \bigg/ \bigg(\frac{|x|}{s}\bigg)^2 \bigg]^{-1} + \bigg[ \frac{V_i(\sqrt t)}{V_i(s)} \bigg/ \bigg(\frac{\sqrt t}{s}\bigg)^2 \bigg]^{-1} \,  \frac{V_i(\sqrt t)}{V_i(y, \sqrt t)} \bigg\}
\nonumber \\
&&\lesssim  \frac{1}{t} \, \frac{s^{2}}{V_i(s)}, \qquad \forall x \in E_i  \setminus F_i^{(s)}, \ y \in F_i^{(s)}, \ t \ge s^2, \ s \ge 1,
\end{eqnarray}}
where we have used in the last inequality \eqref{vg1}, the doubling property and the fact that
\begin{gather*}
d(o_i, y) \le d(y, E_0) + \mbox{diam}\, E_0 \le 2 s + 1, \\
B_i(o_i, \sqrt{t}) \subset B_i(y, \sqrt{t} + d(o_i, y)) \subset B_i(y, 4 \sqrt{t}), \qquad \forall t \ge s^2 \ge 1.
\end{gather*}

Similarly, for $1 \le i \neq  j \le \ell$, \eqref{HEb} implies that
\begin{eqnarray} \label{Pe2}
h_t(x,y) &&\lesssim \frac{1}{t} \, \frac{|x|^2}{V_j(|x|)} + \frac{1}{V_{j}(\sqrt t)} = \frac{1}{t} \, \frac{s^{2}}{V_j(s)} \, \bigg\{ \bigg[ \frac{V_j(|x|)}{V_j(s)} \bigg/ \bigg(\frac{|x|}{s}\bigg)^2 \bigg]^{-1} + \bigg[ \frac{V_j(\sqrt t)}{V_j(s)} \bigg/ \bigg(\frac{\sqrt t}{s}\bigg)^2 \bigg]^{-1} \bigg\} \nonumber \\
&&\lesssim \frac{1}{t} \, \frac{s^{2}}{V_j(s)}, \qquad \forall x \in E_j  \setminus F_j^{(s)}, \ y \in F_i^{(r)}, \ t \ge s^2, \ s, r \ge 1.
\end{eqnarray}

Hence, we obtain that
$$\bigg\| e^{-t\L} \bigg\|_{L^1(F_i^{(R_i)}) \to L^\infty(E_j \setminus F_j^{(R_j)})} \lesssim
\frac{1}{t} \, \frac{R_j^{2}}{V_j(R_j)}, \qquad \forall t \ge R_j^2, \ \gamma \gg 1.$$

Next, recall that $\| e^{- t \L} \|_{L^1 \to L^1} \le 1$.  Hence, the Riesz-Thorin interpolation theorem implies that
$$\bigg\| e^{-t\L} \bigg\|_{L^1(F_i^{(R_i)}) \to L^2(E_j \setminus F_j^{(R_j)})} \lesssim \frac{1}{\sqrt{t}} \, \sqrt{\frac{R_j^{2}}{V_j(R_j)}}, \qquad \forall t \ge R_j^2, \  \gamma \gg 1.$$
By the fact that $\mu(F_i^{(R_i)}) \sim V_i(R_i) \sim V_j(R_j)$
(cf. \eqref{brN}),
then it follows from the H\"older inequality that
$$\bigg\| e^{-t\L} \bigg\|_{L^p(F_i^{(R_i)}) \to L^\infty(E_j \setminus F_j^{(R_j)})} \lesssim
\frac{1}{t} \, \frac{R_j^{2}}{V_j(R_j)}\mu(F_i^{(R_i)})^{1-\frac 1p}\sim \frac{R_j^{2}}{t}\mu(F_i^{(R_i)})^{-\frac 1p},  $$
and
$$\bigg\| e^{-t\L} \bigg\|_{L^p(F_i^{(R_i)}) \to L^2(E_j \setminus F_j^{(R_j)})} \lesssim  \frac{1}{\sqrt{t}} \,
\sqrt{\frac{R_j^{2}}{V_j(R_j)}}  \, \mu(F_i^{(R_i)})^{1-\frac 1p} \sim \sqrt{\frac{R_j^{2}}{t} } \,  \mu(F_i^{(R_i)})^{\frac 12-\frac 1p}.$$
This completes the proof.
\end{proof}

\begin{rem}\rm
From the proof above,  via the doubling property, we see that the assumption that $t \ge R_j^2$ can be relaxed as $t \ge (c R_j)^2$ for some given positive constant $c$.
The same is valid for Proposition \ref{time-heat-map} below.
\end{rem}

\subsection{Mapping properties of $t\L e^{-t\L} = - t \, \frac{d}{d t} e^{-t\L}$ with large time}
\hskip\parindent
For this purpose, let us begin with the follows estimates for the time derivative of the heat kernel in general setting obtained by Davies \cite[Theorem~4]{davies97}.

\begin{thm}[Davies]\label{davies}
Suppose that $\delta\in (0,1)$, $\epsilon\in (0,\frac 18)$, $x,y\in M$ and $t>0$.
Let $a,b,c$ be positive constants such that $c\in (0,1]$, and that
$$h_{(1-\delta) t}(x,x)\le a, \quad h_{(1-\delta)t}(y,y)\le b, \quad |h_{s}(x,y)|\le \sqrt{ab}c$$
for all $s\in ((1-\delta)t,(1+\delta)t)$. Then for any $m\in\cn$, it holds
$$\left|\frac{\partial^m}{\partial t^m}h_t(x,y)\right|\le \frac{m!}{(\epsilon\delta t)^m} \sqrt{ab}c^{1-3\epsilon}. $$
\end{thm}

The following estimates will play a key role in the study of mapping properties for $t\L e^{-t\L}$ with $t \ge 1$:

\begin{lem} \label{TDL}
Let $\epsilon\in (0,1/8)$. Under the assumptions of Theorem \ref{main-result}.   Then:

(i) We have that for $1\le i\le \ell$
\begin{eqnarray} \label{td1}
\bigg| t \, \partial_t h_t(x,y) \bigg| &&\lesssim_\epsilon \frac{1}{t}  \frac{s^2}{V_i(s)} \left(1+ \left(\frac{V_i(s)}{s^2} \right)^{\frac{3}{2} \epsilon}
\left( \frac{|y|^2}{V_i(|y|)} \right)^{\frac{3}{2} \epsilon} \right) \\
&&\lesssim_\epsilon
V_i(s)^{\frac{3}{2} \epsilon - 1} \left( \frac{1}{V_i(|y|)} \right)^{\frac{3}{2} \epsilon},
\  \forall  t \ge s^2 \ge 1,   x \in E_i \setminus F_i^{(s)},   y \in F_i^{(s)}. \nonumber
\end{eqnarray}

(ii) It holds that for all $1 \le i \neq j \le \ell$
\begin{eqnarray} \label{tdN}
\bigg| t \, \partial_t h_t(x,y) \bigg| \lesssim_\epsilon  V_j(s)^{\frac{3}{2} \epsilon - 1} \left( \frac{1}{V_i(|y|)} \right)^{\frac{3}{2} \epsilon}, \qquad \forall t \ge s^2, x \in E_j \setminus F_j^{(s)}, y \in F_i^{(r)}, r, s \ge 1.
\end{eqnarray}
\end{lem}

\begin{proof}
Let us begin with the proof of \eqref{td1}. Recall \eqref{HEa} and \eqref{com-metric}. From the proof of \eqref{Pe1}, it is not hard to notice that
\begin{align} \label{E1}
h_t(x, y) \lesssim \frac{1}{t} \, \frac{s^2}{V_i(s)} \, e^{- c \frac{d(x, y)^2}{t}} := A(x, y, t), \qquad \forall x \in E_i  \setminus F_i^{(s)}, \ y \in F_i^{(s)}, \ t \ge s^2/2  \ge 1/2.
\end{align}

Similarly, we have for all $t \ge s^2/2  \ge 1/2$, $x \in E_i\setminus F_i^{(s)}$ and $y \in F_i^{(s)}$
\begin{gather} \label{E2}
h_t(x, x) \lesssim \frac{1}{t} \, \frac{|x|^2}{V_i(|x|)} +  \frac{1}{V_i(x, \sqrt{t})} \lesssim \frac{1}{t} \, \frac{s^2}{V_i(s)} +  \frac{1}{V_i(x, \sqrt{t})} := a(x, t), \\
h_t(y, y) \lesssim \frac{1}{t} \, \left\{ \frac{|y|^2}{V_i(|y|)}  +  \frac{s^2}{V_i(s)} \right\} = \frac{1}{t} \, \frac{s^2}{V_i(s)} \, \left\{ 1 + \frac{V_i(s)}{s^2} \,  \frac{|y|^2}{V_i(|y|)}  \right\} := b(y, t).
  \label{E3}
\end{gather}

From \eqref{E1}-\eqref{E3},  Theorem \ref{davies} implies that
\begin{align*}
\bigg| t \, \frac{\partial}{\partial t} h_t(x, y) \bigg| &\lesssim_{\epsilon} A(x, y, t)^{1 - 3 \epsilon} \left( a(x, t) \, b(y, t) \right)^{\frac{3}{2} \epsilon} \\
&\le \left( \frac{1}{t} \, \frac{s^2}{V_i(s)} \right)^{1 - 3 \epsilon} \, e^{- c' \frac{d(x, y)^2}{t}} \,  \left( a(x, t) \, b(y, t) \right)^{\frac{3}{2} \epsilon}.
\end{align*}
Next we decompose $e^{- c' \frac{d(x, y)^2}{t}}$ as $e^{- c \, \frac{d(x, y)^2}{t}} \, e^{- c \, \frac{d(x, y)^2}{t}}$, and write
\begin{align*}
e^{- c \, \frac{d(x, y)^2}{t}} \, a(x, t) \le \frac{1}{t} \, \frac{s^2}{V_i(s)} +  \frac{1}{V_i(x, \sqrt{t})} \, e^{- c \, \frac{d(x, y)^2}{t}} \lesssim  \frac{1}{t} \, \frac{s^2}{V_i(s)} +  \frac{1}{V_i(y, \sqrt{t})} \lesssim \frac{1}{t} \, \frac{s^2}{V_i(s)},
\end{align*}
where we have used { in the second inequality the standard trick of doubling property (cf. the proof of \eqref{BE1}),} and \eqref{Pe1} again in the last inequality. In conclusion, by the trivial inequality $(1 + a)^b \lesssim 1 + a^b$ for all $a > 0$ and $0 < b \le 1$, we get that
\begin{align*}
\bigg| t \, \frac{\partial}{\partial t} h_t(x, y) \bigg| &\lesssim_{\epsilon} \frac{1}{t} \, \frac{s^2}{V_i(s)} \, e^{- c \frac{d(x, y)^2}{t}} \,  \left(1+ \left(\frac{V_i(s)}{s^2} \right)^{\frac{3}{2} \epsilon}
\left( \frac{|y|^2}{V_i(|y|)} \right)^{\frac{3}{2} \epsilon} \right),
\end{align*}
which is clearly majorized by constant times
$V_i(s)^{- 1} \left(V_i(s)/V_i(|y|) \right)^{\frac{3}{2} \epsilon}$ under our assumptions, by the definition of $F_i^{(r)}$ (cf. \eqref{na}) and the doubling property.

Now we turn to the proof of \eqref{tdN}.  In a very similar way,  we have that for all $t \ge s^2/2 \ge 1/2$, $x \in E_j \setminus F_j^{(s)}$ and $y \in F_i^{(r)}$  with $1 \le i \neq j \le \ell$,
\begin{gather*}
h_t(x, y)  \lesssim \frac{1}{t} \, \frac{s^{2}}{V_j(s)}  \, e^{- c \frac{|x|^2+|y|^2}{t}} :=A_*(x, y, t), \qquad \mbox{from the proof of \eqref{Pe2},} \\
h_t(x, x) \lesssim \frac{1}{t} \, \frac{s^2}{V_j(s)} +  \frac{1}{V_j(x, \sqrt{t})} := a_*(x, t), \\
h_t(y, y) \lesssim \frac{1}{t} \, \frac{|y|^2}{V_i(|y|)}  +  \frac{1}{V_i(y, \sqrt{t})}  := b_*(y, t).
\end{gather*}
Then, it deduces from Theorem \ref{davies} that
\begin{align} \label{F4}
\bigg| t \, \frac{\partial}{\partial t} h_t(x, y) \bigg| \lesssim_{\epsilon}
\left( \frac{1}{t} \, \frac{s^2}{V_j(s)} \right)^{1 - 3 \epsilon} \, e^{- c \frac{|x|^2 + |y|^2}{t}} \,  \left( e^{- c \frac{|x|^2}{t}} a_*(x, t) \, e^{- c \frac{|y|^2}{t}} b_*(y, t) \right)^{\frac{3}{2} \epsilon}.
\end{align}

Now by \eqref{vg1}, the standard trick of doubling property implies that (e.g. in the proof of Corollary \ref{uE})
\begin{align*}
e^{- c \frac{|x|^2}{t}} a_*(x, t) \lesssim \frac{1}{t} \, \frac{s^2}{V_j(s)}, \qquad e^{- c \frac{|y|^2}{t}} b_*(y, t) \lesssim \frac{1}{t} \, \frac{|y|^2}{V_i(|y|)}.
\end{align*}
Combining this with \eqref{F4}, we obtain that
\begin{align*}
\bigg| t \, \frac{\partial}{\partial t} h_t(x, y) \bigg| \lesssim_{\epsilon}  \left( \frac{1}{t} \, \frac{s^2}{V_j(s)} \right)^{1 - \frac{3}{2} \epsilon} \,  \left( \frac{1}{t} \, \frac{|y|^2}{V_i(|y|)} \right)^{\frac{3}{2} \epsilon} \, e^{- c \frac{|y|^2}{t}},
\end{align*}
then by the trivial inequality $e^{- c \frac{|y|^2}{t}}  \lesssim_{c, \epsilon} \left( \frac{|y|^2}{t} \right)^{- \frac{3}{2} \epsilon}$, the fact that $t \ge s^2$ implies the desired result.
\end{proof}

\begin{lem}
Suppose that $0 < \frac{3}{2} \epsilon < \frac{1}{p'}$. Then we have for all $1 \le i \le \ell$ and $s \ge 1$
\begin{align}  \label{VeN}
I(i, \epsilon, p', s) := \left\{ \int_{\{y \in E_i; \, d(y, E_0) \le 2 s \}} \left(  \frac{1}{V_i(|y|)}  \right)^{\frac{3}{2} \epsilon p'} \, d\mu(y) \right\}^{\frac{1}{p'}} \lesssim_{\epsilon, p'} V_i(s)^{\frac{1}{p'} - \frac{3}{2} \epsilon}.
\end{align}
\end{lem}

\begin{proof}
We split the domain of integration into the (not necessarily disjoint) sets
\[
\{y \in E_i; \, d(y, E_0) \le 2 \}, \qquad \{y \in E_i; \, 2^{- k} s < d(y, E_0) \le 2^{-k + 1} s \}  \quad \mbox{with } \ 0 \le k \le [\log_2 s],
\]
where $[a]$ denotes the integer part of $a \in \rr$. Hence, by \eqref{vcp}, we can write
\[
I(i, \epsilon, p', s) \lesssim 1 + \sum_{k = 0}^{[\log_2 s]} \left(  \frac{1}{V_i(2^{-k} s)}  \right)^{\frac{3}{2} \epsilon} \, V_i(2^{-k + 1} s)^{\frac{1}{p'}}.
\]

From the doubling property and the facts that $\frac{3}{2} \epsilon  < \frac{1}{p'}$, we get that
\begin{eqnarray*}
I(i, \epsilon, p', s) \lesssim 1 + \sum_{k = 0}^{[\log_2 s]} V_i(2^{-k} s)^{\frac{1}{p'} - \frac{3}{2} \epsilon} &&\le 2 \, V_i(s)^{\frac{1}{p'} - \frac{3}{2} \epsilon} \sum_{k = 0}^{[\log_2 s]} \left(  \frac{V_i(s)}{V_i(2^{-k} s)} \right)^{\frac{3}{2} \epsilon - \frac{1}{p'}} \\
&&\lesssim V_i(s)^{\frac{1}{p'} - \frac{3}{2} \epsilon} \sum_{k = 0}^{+\infty} 2^{- n_i k ( \frac{1}{p'} - \frac{3}{2} \epsilon)}\\
&&\lesssim_{\epsilon, p'}  V_i(s)^{\frac{1}{p'} - \frac{3}{2} \epsilon},
\end{eqnarray*}
where the penultimate inequality follows from \eqref{vg1}.
\end{proof}

Finally one can state
\begin{prop}\label{time-heat-map}
Under the assumptions of Theorem \ref{main-result}, for each $p\in (1,2)$, it holds  for all $1 \le i, j \le \ell$ that
\begin{align} \label{Ntd}
\left\|t\L e^{-t\L}\right\|_{L^p(F_i^{(R_i)}) \to L^2(E_j \setminus F_j^{(R_j)})} \lesssim_p
\mu(F_i^{(R_i)})^{\frac 12-\frac 1p}, \qquad \forall t \ge R_j^2=R_j(\gamma)^2, \  \gamma \gg 1.
\end{align}
\end{prop}
\begin{proof}
It follows from  Lemma \ref{TDL} that
\begin{align*}
\left\| t \L e^{-t\L}(f \chi_{F_i^{(R_i)}}) \right\|_{L^\infty(E_j \setminus F_j^{(R_j)})} \lesssim_{\epsilon} V_j(R_j)^{\frac{3}{2} \epsilon - 1} \int_{F_i^{(R_i)}} |f(y)| \, \left( \frac{1}{V_i(|y|)} \right)^{\frac{3}{2} \epsilon} \, d\mu(y).
\end{align*}
Then by the  H\"older inequality, we obtain from \eqref{VeN} that
$$\left\| t\L e^{-t\L}(f \chi_{F_i^{(R_i)}}) \right\|_{L^\infty(E_i \setminus F_i^{(R_i)})}  \lesssim_{\epsilon} V_j(R_j)^{\frac{3}{2} \epsilon - 1} \, V_i(R_i)^{\frac{1}{p'} - \frac{3}{2} \epsilon} \, \| f \|_{L^p(F_i^{(R_i)})} \lesssim_{\epsilon} \frac {\|f\|_{L^p(F_i^{(R_i)})}}{\mu(F_i^{(R_i)})^{1/p}},$$
since $\mu(F_i^{(R_i)}) \sim V_i(R_i) \sim V_j(R_j)$  (cf. \eqref{brN}).

In conclusion, for any $1< p < 2$, by suitably choosing $\epsilon$, we have that
\begin{eqnarray*}
&&\left\| t\L e^{-t\L} \right\|_{L^p(F_i^{(R_i)}) \to L^\infty(E_j \setminus F_j^{(R_j)})} \lesssim_p \frac{1}{\mu(F_i^{(R_i)})^{1/p}}.
\end{eqnarray*}
On the other hand, the classical Littlewood-Paley-Stein theory says that
\begin{eqnarray*}
&&\left\| t \L e^{-t\L} \right\|_{L^p(F_i^{(R_i)}) \to L^p(E_j \setminus F_j^{(R_j)})}
\lesssim_p 1, \qquad \forall\, 1 < p < +\infty.
\end{eqnarray*}

Hence the H\"older inequality implies that
\begin{eqnarray*}
\left\| t \L e^{-t\L} \right\|_{L^p(F_i^{(R_i)}) \to L^2(E_j \setminus F_j^{(R_j)})} \lesssim_p \mu(F_i^{(R_i)})^{\frac 12-\frac 1p},
\end{eqnarray*}
which completes the proof.
\end{proof}

\subsection{{$L^1$ estimate of the space derivative of the heat kernel}}
\hskip\parindent
The following result can be considered as a counterpart of the main result in \cite[\S~2.3]{cd99}, in other words, weighted estimates for the space derivative of the heat kernel (cf. \cite[Lemma~2.4]{cd99}). However, we will provide a direct proof by means of the main results in \cite{G95, G97}. {And another proof similar to  \cite[Lemma~2.4]{cd99} can be found in Subsection 4.5 below.

\begin{prop}\label{est-integral-diagonal}
Under the assumptions of Theorem \ref{main-result}. Let $\beta \ge 16$ and $1 \le i \le \ell$,  it holds  that:
\begin{eqnarray*}
\int_{E_i \setminus (F_i^{(s)} \cup B(y, r))}|\nabla_x h_t(x, y)| \,d\mu(x) \lesssim_{\beta} \frac{1}{\sqrt t} \, e^{- c \frac{r^2}{t}}, \qquad \forall t \le \beta s^2,  \ y \in E_i \setminus F_i^{(s)}, \ s \ge 1,  \ r > 0.
\end{eqnarray*}
\end{prop}
\begin{proof}
First by the H\"older inequality, we can write that
\begin{eqnarray*}
\int_{E_i \setminus (F_i^{(s)} \cup B(y,r))} | \nabla_x h_t(x, y)|  \, d\mu(x) &&\le \left(\int_M |\nabla_x h_t(x, y)|^2 \, e^{\frac{d(x, y)^2}{3 t}} \, d\mu(x) \right)^{\frac{1}{2}} \left(\int_{E_i \setminus B(y, r)} e^{-\frac{d(x, y)^2}{3 t}} \, d\mu(x) \right)^{\frac{1}{2}} \\
&&=: E_1(y, t)^{\frac{1}{2}} \, G(y, t, r)^{\frac{1}{2}}.
\end{eqnarray*}

Next, the doubling property and a standard method of decomposition in annuli imply that (e.g. \cite[Lemma~2.1]{cd99})
\begin{align} \label{WB1}
G(y, t, r) \lesssim e^{-\frac{r^2}{6 t}} \, V_i(y, \sqrt{t}), \qquad \forall t, r > 0, \  y \in E_i.
\end{align}

Moreover Corollary \ref{uE} says that
\[
h_t(y, y) \lesssim_{\beta} \frac{1}{V_i(y, \sqrt{t})}, \quad \forall 0 < t \le \beta s^2, \  y \in E_i \setminus F_i^{(s)}, \ s \ge 1.
\]
Using the doubling property, it is clear that the function $h(t) := h(y, t) :=  V_i(y, \sqrt{t})$ is regular on $(0, \, \beta R_i^2)$ in the sense of \cite[p.~37]{G97} and $h^{(1)}(t) := \int_0^t h(s) \, ds \sim t \, h(t)$ (e.g. \cite[Lemma~3.1]{G95}). We then apply \cite[Theorem~3.1]{G97} and \cite[Theorem~1.1]{G95}, and get that
\begin{align} \label{WB2}
E_1(y, t) \lesssim_{\beta} t^{-1} \frac{1}{V_i(y, \sqrt{t})}, \qquad \forall 0 < t \le \beta s^2, \  y \in E_i \setminus F_i^{(s)}, \ s \ge 1.
\end{align}
Combining this with \eqref{WB1}, we obtain the desired estimates.
\end{proof}

}

\section{Proof of Theorem \ref{main-result}}

\subsection{Outline of the proof}
\hskip\parindent  Let us outline the main approach for convenience of the reader.
The heat kernel estimates { (cf. Theorem \ref{hke})} established in  Grigor'yan and Saloff-Coste \cite{gri-sal09} is a key tool for the proof.  Since we showed in Lemma \ref{local-part} that the operator
$$\frac{1}{\sqrt{\pi}} \int_0^1\nabla e^{-s\L}\frac{\,ds}{\sqrt s}$$
is bounded on $L^p(M)$ for $1<p<2$, we only need to show that
$$T:=\frac{1}{\sqrt{\pi}} \int_1^\infty\nabla e^{-s\L}\frac{\,ds}{\sqrt s}$$
is bounded on $L^p(M)$ for $1<p<2$. To this end, we will show that $T$ is weakly $(p,p)$ bounded for any $1<p<p_0$, here $p_0$ is defined as
\begin{equation}\label{set-index-1}
p_0:=\min\left\{2,\,\frac{N_\infty}
{N_\infty-2}\right\},
\end{equation}
where $N_\infty=\max_i\left\{N_i:\,i=1,\cdots,\ell\right\}$ is defined in \eqref{set-index}.

{\bf Step 1: Reduction}.
For $1<p<p_0$ and any non-trivial $f\in L^p(M)$, from the non-collapsing assumption,
it is not hard to see that for any $\lambda$ such that $\|f\|_{L^p(M)}\lesssim \lambda$,
$$\mu\left(\left\{x:\, |Tf|>\lambda\right\}\right)\lesssim \frac{\|f\|_p^p}{\lambda^p};$$
see Subsection 4.2. It is then enough to prove the same estimate holds for
the case $0<\lambda \ll \|f\|_{L^p(M)}$. To this end, in each end $E_i$, we choose a set $F_i^{(R_i)}$
around $E_0$ such that $\mu(F_i^{(R_i)})\sim \|f\|^p_{L^p(M)}/\lambda^p$; see Figure 2.
It is then enough to estimate that
\begin{eqnarray*}
&&\mu\left(\left\{x\in \bigcup_{j=1}^\ell \left( E_j\setminus F_j^{(R_j)} \right): |Tf|>(2\ell+1)\lambda\right\}\right)\\
&&\quad\le \mu\left(\left\{x\in \bigcup_{j=1}^\ell \left( E_j\setminus F_j^{(R_j)} \right): |T(f\chi_{E_0})|>\lambda\right\}\right)+\sum_{i=1}^\ell
\mu\left(\left\{x\in \bigcup_{j=1}^\ell \left( E_j\setminus F_j^{(R_j)} \right): |T(f\chi_{F_i^{(R_i)}})|>\lambda\right\}\right)\\
&&\quad\quad+\sum_{i=1}^\ell
\mu\left(\left\{x\in \bigcup_{j=1}^\ell \left( E_j\setminus F_j^{(R_j)} \right): |T(f\chi_{E_i\setminus F_i^{(R_i)}})|>\lambda\right\}\right).
\end{eqnarray*}

{\bf Step 2: Estimate of the center part}.
The part $f\chi_{E_0}$ is obvious, see Subsection 4.3.

{\bf Step 3: Estimate of the part away from the center}.
For the part $f\chi_{E_i\setminus F_i^{(R_i)}}$,
 we shall incorporate some ideas from \cite{cd99} and use the Calder\'on-Zygmund decomposition to decompose
 $f\chi_{E_i\setminus F_i^{(R_i)}}$ as $g_i+\sum_k b_{ik}$, where $g_i$ and $b_{ik}$ are supported on $E_i$,  and consider
$$
\mu\left(\left\{x\in \bigcup_{j\neq i} \left( E_j\setminus F_j^{(R_j)} \right): |T(f\chi_{E_i\setminus F_i^{(R_i)}})|>\lambda\right\}\right),
\, \& \, \mu\left(\left\{x\in E_i\setminus F_i^{(R_i)}:|T(f\chi_{E_i\setminus F_i^{(R_i)}})|>\lambda\right\}\right).$$
Since we are dealing with small $\lambda$, $\lambda \ll \|f\|_{L^p(M)}$, the supporting balls of
$b_{ik}$ are small enough, the $L^p$-Davies-Gaffney estimate of the operators $\nabla e^{-t\L}$ and $\L e^{-t\L}$ established in Section 2 is sufficient to yield the estimate on the off-diagonal part, i.e., the first term. For the diagonal part (the second term), we shall further decompose the operator $T$ into the parts above or below the diagonal, and deduce the required estimate by using the $L^p$-boundedness (Davies-Gaffney estimate) of  $\nabla e^{-t\L}$ and $\L e^{-t\L}$,  and  heat kernel estimate in Subsection 3.3.

{\bf Step 4: Estimate of the part around the center}.
For the part  $f\chi_{F_i^{(R_i)}}$, we shall use the boundedness of the operators $e^{-t\L}$
and $\L e^{-t\L}$ from $L^p(F_i^{(R_i)})$ to $L^2(E_j\setminus F_j^{(R_j)})$ in Section 3,
and the $L^p$ Davies-Gaffney estimate from Section 2.

{\bf Step 5: Completion of the proof}. By combining previously obtained results,
we show that $\nabla \L^{-1/2}$ is weakly $(p,p)$ bounded for each $1<p<p_0$.
An application of the Marcinkowicz interpolation theorem gives the desired result.

\subsection{Reduction}
\hskip\parindent
Let us begin with the following simple observation:

\begin{lem}\label{est-large}
Let $1 < p < 2$ and $0 < \beta_0 \le 1$. It holds that
\begin{eqnarray}\label{est-large-lambda}
\mu\left(\left\{x:\, |Tf|>\lambda\right\}\right)  \lesssim_{\beta_0} \frac{\|f\|_p^p}{\lambda^p}, \qquad \forall \lambda \ge \beta_0 \|f\|_{p}, \  f \in L^p.
\end{eqnarray}
\end{lem}

\begin{proof}
Using Chebyshev's inequality and Minkowski's integral  inequality, we have that
\begin{eqnarray*}
\mu\left(\left\{x:\, |Tf|>\lambda\right\}\right)&&\lesssim \frac{1}{\lambda^2}\left\|\int_1^\infty\nabla e^{-t\L}f\frac{\,dt}{\sqrt t}\right\|^2_2 \lesssim  \frac{1}{\lambda^2}\left[ \int_1^\infty\left\|\nabla e^{-t\L}f\right\|_2 \frac{\,dt}{\sqrt t}\right]^2.
\end{eqnarray*}
Next, the semigroup property and Proposition \ref{mapping-gradient-heat} imply that
\[
\left\|\nabla e^{-t\L}f\right\|_2 = \left\|\nabla e^{-\frac{t}{2}\L} \circ e^{-\frac{t}{2}\L} f\right\|_2 \lesssim \frac{1}{\sqrt{t}} \left\| e^{-\frac{t}{2}\L} f\right \|_2 \lesssim t^{\delta (\frac{1}{2} - \frac{1}{p})} \frac{1}{\sqrt{t}} \| f \|_p
\]
because of the ultracontractivity of the heat semigroup (cf. \eqref{heat-operator-norm}).
In conclusion,
\begin{eqnarray*}
\mu\left(\left\{x:\, |Tf|>\lambda\right\}\right) \lesssim \frac{\|f\|_p^2}{\lambda^2} \lesssim_{\beta_0} \frac{\|f\|_p^p}{\lambda^p}, \quad \forall \lambda \ge \beta_0 \| f \|_p,
\end{eqnarray*}
which finishes the proof.
\end{proof}

It is then enough to prove \eqref{est-large-lambda} holds for  the case $0<\lambda \ll \|f\|_{p}$ with $1 < p < p_0$ (see \eqref{set-index-1} for $p_0$).  In such case,
we decompose $M$ as
$$M=E_0\cup (\cup_{i=1}^\ell F_i^{(R_i)}) \cup {\Big[ \cup_{j=1}^\ell (E_j\setminus F_j^{(R_j)}) \Big]}$$
where
\begin{align} \label{fi}
F_i^{(R_i)} := \{x\in E_i:\,\dist(x,E_0)\le 2R_i\} \  \mbox{with $R_i \gg 1$ such that} \
\mu(F_i^{(R_i)}) = 100 \, \ell \, \frac{\|f\|_p^p}{\lambda^p};
\end{align}
see Figure 2. By the convention $\mu(E_0) = 1$ (cf. \eqref{c1}), it holds then
\begin{eqnarray}\label{est-small-lambda-1}
\mu{\Big( E_0\cup (\cup_{i=1}^\ell F_i^{(R_i)}) \Big)} = 1+\sum_{i=1}^\ell \mu(F_i^{(R_i)})\sim \frac{\|f\|_p^p}{\lambda^p}.
\end{eqnarray}

Therefore, to prove that \eqref{est-large-lambda} holds for such $\lambda$, it is enough to prove that
\begin{eqnarray}\label{est-small-lambda-2}
&&\mu\left(\left\{ x\in \bigcup_{j=1}^\ell \left( E_j\setminus F_j^{(R_j)} \right):|Tf|>\lambda \right\}\right)\lesssim  \frac{\|f\|_p^p}{\lambda^p}.
\end{eqnarray}
For the simplicity of notation, we replace $\lambda$ by $(2\ell+1)\lambda$ and write
\begin{eqnarray*}
&&\mu\left(\left\{ x\in \bigcup_{j=1}^\ell \left( E_j\setminus F_j^{(R_j)} \right):|Tf|>(2\ell+1)\lambda \right\}\right)\\
&&\quad\le \mu\left(\left\{ x\in \bigcup_{j=1}^\ell \left( E_j\setminus F_j^{(R_j)} \right):|T(f\chi_{E_0})|>\lambda \right\}\right)+\sum_{i=1}^\ell
\mu\left(\left\{ x\in \bigcup_{j=1}^\ell \left( E_j\setminus F_j^{(R_j)} \right):|T(f\chi_{F_i^{(R_i)}})|>\lambda \right\}\right)\\
&&\quad\quad+\sum_{i=1}^\ell
\mu\left(\left\{ x\in \bigcup_{j=1}^\ell \left( E_j\setminus F_j^{(R_j)} \right):|T(f\chi_{E_i\setminus F_i^{(R_i)}})|>\lambda \right\}\right)
=: \ncc + \ca + \ncw.
\end{eqnarray*}
{  The same trick will be used repeatedly.}
We shall call the first term as the center part, the second term as the part around the center, and
the last term as the part away from the center. We shall treat them in the following subsections.

\subsection{Estimate on the center}\hskip\parindent
First we establish the following:
\begin{lem}\label{est-center}
Let $1<p<2$. It holds that
\begin{eqnarray*}
 \ncc = \mu\left(\left\{ x\in \bigcup_{j=1}^\ell \left( E_j\setminus F_j^{(R_j)} \right): |T(f\chi_{E_0})|>\lambda \right\}\right)
\lesssim \frac{\|f\|_p^p}{\lambda^p}, \qquad \forall \lambda \ll \|f\|_{p}, \  f \in L^p.
\end{eqnarray*}
\end{lem}
\begin{proof}
Similar to the proof of \eqref{est-large-lambda}, we have that
\begin{eqnarray*}
\ncc \lesssim \frac{1}{\lambda^p}\left\|T(f\chi_{E_0})\right\|_p^p
&&\lesssim \frac{1}{\lambda^p}\left[\int_1^\infty \left\| \nabla e^{-t\L}(f\chi_{E_0}) \right\|_p\frac{\,dt}{\sqrt t}\right]^p
\\
&&\lesssim \frac{1}{\lambda^p}\left[\int_1^\infty  \left\| e^{-\frac{t}{2}  \L}(f\chi_{E_0}) \right\|_p\frac{\,dt}{t}\right]^p\\
&&\lesssim \frac{1}{\lambda^p}\left[\int_1^\infty t^{\delta(\frac 1p-1)}\|f\chi_{E_0})\|_1\frac{\,dt}{t}\right]^p\\
&&\lesssim \frac{\|f\chi_{E_0}\|_1^p}{\lambda^p}\lesssim \frac{\|f\|_p^p}{\lambda^p},
\end{eqnarray*}
as desired.
\end{proof}
\subsection{The part away from the center}\hskip\parindent
Recall that $p_0$ is defined by \eqref{set-index-1}.
The main aim of this subsection is to prove that
\begin{prop}\label{est-away-center}
Let $1<p<p_0$. It holds that
$${\ncw = } \sum_{i=1}^\ell
\mu\left(\left\{ x\in \bigcup_{j=1}^\ell \left( E_j\setminus F_j^{(R_j)} \right):|T(f\chi_{E_i\setminus F_i^{(R_i)}})|>\lambda \right\}\right)\lesssim \frac{\|f\|_p^p}{\lambda^p}, \qquad \forall \lambda \ll \|f\|_{p}, \  f \in L^p.$$
\end{prop}
For each part $E_i\setminus F_i^{(R_i)}$, we may run the Calder\'on-Zygmund decomposition (cf. \cite{cw77}) of $f\chi_{E_i\setminus F_i^{(R_i)}}$ on $M_i$ at the height $\delta\lambda$, where $\delta>1$ is large enough to be fixed later; see also \cite{cd99}.  Since $M_i$ is a doubling manifold, we obtain
$$f\chi_{E_i\setminus F_i^{(R_i)}}=\sum_{k}b_{ik}+g_i,$$
where $\supp b_{ik}\subset B_{ik}=B_i(x_{ik},r_{ik})$,
 $\sum_k \chi_{B_{ik}}\lesssim 1$,  $|g_i|\le \delta\lambda$,
\begin{equation}\label{cz-1}
\int_{M_i}|g_i|^p\,d\mu_i \lesssim \int_{M_i}|f\chi_{E_i\setminus F_i^{(R_i)}}|^p\,d\mu_i = \int_{M}|f\chi_{E_i\setminus F_i^{(R_i)}}|^p\,d\mu,
\end{equation}
\begin{equation}\label{cz-2}
\fint_{B_{ik}}|b_{ik}|^p\,d\mu_i\lesssim  (\delta\lambda)^p,
\end{equation}
and
\begin{equation}\label{cz-3}
\sum_{k}\mu(B_{ik})\lesssim \frac{1}{(\delta\lambda)^p}\int_{M_i}|f\chi_{E_i\setminus F_i^{(R_i)}}|^p\,d\mu_i.
\end{equation}
We shall call $g_i$ the good part, and $b_{ik}$ the `bad' part following tradition.

From \eqref{cz-3} we further have
\begin{equation}\label{cz-measure}
\sum_{k}\mu(B_{ik})\lesssim \frac{1}{(\delta\lambda)^p}\int_{M_i}|f {\chi_{E_i \setminus F_i^{(R_i)}} }|^p\,d\mu_i \lesssim \frac{\|f\|_p^p}{(\delta\lambda)^p}\lesssim\frac{\mu(F_i^{(R_i)})}{\delta^p},
\end{equation}
by \eqref{fi}. We claim that, for a large enough $\delta$, it holds that
\begin{equation}\label{isolate-balls}
\dist(B_{ik},E_0)=\dist(B_{ik},\partial E_i)\ge \frac{1}{2}(R_i+r_{ik}) (\gg 1).
\end{equation}
 Indeed, if the radius $r_{ik}<R_i/2$, then the claim holds obviously, since $B_{ik}\cap (E_i\setminus F_i^{(R_i)}) \neq\emptyset$ and $F_i^{(R_i)}=\{x\in E_i:\,\dist(x,E_0)\le 2R_i\}$.   If $r_{ik}\ge R_i/2$ and the claim fails, then
\begin{equation*}
\dist(B_{ik},E_0)<\frac{1}{2}(R_i+r_{ik})<\frac{3}{2}r_{ik}.
\end{equation*}
This together with the convention $\mbox{diam}(E_0)=1$ (cf. \eqref{c1}) and the doubling condition implies that
\begin{equation*}
{ \mu\left( F_i^{(R_i)} \right) \sim} V_i(R_i)\lesssim V_i(o_i,3r_{ik})\lesssim V_i(x_{ik},3r_{ik})\lesssim  \mu(B_{ik}),
\end{equation*}
which contradicts with \eqref{cz-measure} as soon as we choose $\delta$ large enough.
Therefore, \eqref{isolate-balls} holds, which further implies that
\begin{equation}\label{isolate-balls-1}
\dist(B_{ik},E_0)\le d(x_{ik},o_i) \le  \dist(B_{ik},E_0)+\diam(E_0)+r_{ik}\le 3\dist(B_{ik},E_0),
\end{equation}
and
\begin{equation}\label{isolate-balls-2}
\frac{1}{2}(R_i+r_{ik})+r_{ik}\le \dist(B_{ik},E_0)+r_{ik}\le d(x_{ik},o_i).
\end{equation}

Note that \eqref{isolate-balls} implies that $\supp b_{ik}\subset E_i$, which together with the identity
$f\chi_{E_i\setminus F_i^{(R_i)}}=\sum_{k}b_{ik}+g_i$ implies that
$\supp g_i\subset E_i.$
Since $\mu=\mu_i$ on $E_i$, the integrals  of $g_i$ and $b_{ik}$ on $M_i$
are the same as on $M$, i.e.,
$$\int_{M_i}|g_i|^p\,d\mu_i= \int_{M}|g_i|^p\,d\mu, \qquad
\int_{B_{ik}}|b_{ik}|^p\,d\mu_i=\int_{B_{ik}}|b_{ik}|^p\,d\mu.$$

\subsubsection{Estimate of the good part}
\hskip\parindent The $L^2$-boundedness of $T$ (cf. Lemma \ref{local-part} (i))
together with the fact $|g_i|\lesssim \lambda$  implies
\begin{eqnarray}\label{est-g}
\sum_{i=1}^\ell
\mu\left( \left\{ x\in \bigcup_{j=1}^\ell \left( E_j\setminus F_j^{(R_j)} \right) :|T(g_i)|>\lambda \right\} \right)
&&\lesssim \sum_{i=1}^\ell\frac{1}{\lambda^2}\int_{M}|g_i|^2\,d\mu
\lesssim \sum_{i=1}^\ell\frac{1}{\lambda^p}\int_{M_i}|g_i|^p\,d\mu\nonumber\\
&&\lesssim \frac{\|f\|_p^p}{\lambda^p},
\end{eqnarray}
{  because of \eqref{cz-1}.}

\subsubsection{Estimate of the `bad' parts off the diagonal}
\hskip\parindent For each $i$, $1\le i\le \ell$,
\begin{eqnarray}\label{partion-off-diagonal}
&&
\mu\left(\left\{x\in \bigcup_{1 \le j\neq i \le \ell} \left( E_j\setminus F_j^{(R_j)} \right): \left|T\left(\sum_k b_{ik}\right)\right|>\lambda\right\}\right)\nonumber\\
&&\quad\le
\mu\left(\left\{x\in \bigcup_{1 \le j\neq i \le \ell} \left( E_j \setminus F_j^{(R_j)} \right): \left|T\left(\sum_k (1-e^{-r_{ik}^2\L}) b_{ik}\right)\right|>\lambda/2\right\}\right)\nonumber\\
&&\quad\quad+
\mu\left(\left\{x\in \bigcup_{1 \le j\neq i \le \ell} \left( E_j\setminus F_j^{(R_j)} \right): \left|T\left(\sum_k e^{-r_{ik}^2\L}b_{ik}\right)\right|>\lambda/2\right\}\right)
\nonumber\\
&&\quad=: \cg_1 + \cg_2.
\end{eqnarray}
For the term $\cg_2$, we have
\begin{prop}\label{est-away-center-offdiagonal-1}
Let $1<p<2$. It holds that
\begin{eqnarray*}
\cg_2 =  \mu\left(\left\{x\in \bigcup_{ 1\le j\neq i\le \ell} \left( E_j\setminus F_j^{(R_j)} \right): \left| T\left(\sum_k e^{-r_{ik}^2\L}b_{ik}\right) \right| > \lambda/2 \right\}\right)\lesssim \frac{\|f\|_p^p}{\lambda^p}, \qquad \forall \lambda \ll \|f\|_{p}, \  f \in L^p.
\end{eqnarray*}
\end{prop}
\begin{proof}
The $L^2$-boundedness of $T$ implies that
\begin{eqnarray}\label{est-b-off}
&&  \cg_2 \lesssim \frac{1}{\lambda^2}\int_M \left|\sum_k e^{-r_{ik}^2\L}b_{ik}\right|^2\,d\mu\nonumber\\
&&\quad= \lambda^{-2} \, \int_{E_i} \left|\sum_k e^{-r_{ik}^2\L}b_{ik} \right|^2 \, d\mu + \lambda^{-2} \, \sum_{0 \le j \neq i \le \ell} \int_{E_j} \left|\sum_k e^{-r_{ik}^2\L}b_{ik}\right|^2 \, d\mu  =: \cg_{21} + \cg_{22}.
\end{eqnarray}

Let us begin with $\cg_{21}$. By \eqref{isolate-balls}, it follows from Corollary \ref{uE} that
\begin{eqnarray*}
h_{r_{ik}^2}(x, y) \lesssim  \frac{1}{V_i(x,r_{ik})}e^{- c \frac{d(x, y)^2}{t}}, \qquad \forall  x \in E_i,  \  y\in B_{ik}\subset E_i.
\end{eqnarray*}

Let $\mathcal{M}_i$ denote the centered Hardy-Littlewood Maximal function on $M_i$.
Then the same arguments as Coulhon-Duong \cite[pp. 1158-1160]{cd99} gives for any $\|\psi\|_2\le 1$ that
\begin{eqnarray*}
\int_{E_i} \left( \sum_k e^{-r_{ik}^2\L}b_{ik} \right) \, \psi\,d\mu&&\lesssim  \lambda\int_{E_i} \left( \sum_{k}\chi_{B_{ik}} \right) \, \mathcal{M}_i\psi\,d\mu\lesssim
\lambda\left\|\sum_{k}\chi_{B_{ik}}\right\|_{L^2(M_i)}\|\mathcal{M}_i\psi\|_{L^2(M_i)}\\
&&\lesssim
\lambda\left(\sum_{k}\mu(B_{ik})\right)^{1/2}.
\end{eqnarray*}
Combining this with \eqref{cz-3}, we obtain that
$$\cg_{21} \lesssim \sum_{k}\mu(B_{ik}) \lesssim \frac{\| f \|_p^p}{\lambda^p}.$$

Consequently, what remains is to prove that
$$\cg_{22} = \lambda^{-2} \, \sum_{0 \le j \neq i \le \ell}\int_{E_j}\left|\sum_k e^{-r_{ik}^2\L}b_{ik}\right|^2 \, d\mu \lesssim \frac{\| f \|_p^p}{\lambda^p}, \quad \forall \lambda \ll \| f \|_p, \ f \in L^p, \  1 \le i \le \ell.$$

Next, we claim that
\begin{align} \label{nSGn}
\left\|e^{-r_{ik}^2\L}\right\|_{L^1(M \setminus E_i) \to L^\infty(B_{ik})} \lesssim \frac{1}{V_i(R_i)}.
\end{align}
Let $u \in  B_{ik} \subset E_i$ ($1 \le i \le \ell$) and $v \not\in  E_i$.
First consider the case where $r_{ik}\ge 1$. Using \eqref{h} and the fact that $V_j(s) \ge V_0(s) \gtrsim s^2$ for any $s \ge 1$ and $0 \le j \le \ell$ (cf. Lemma \ref{vg} (i)), Theorem \ref{hke} (iii) implies that
\begin{eqnarray*}
h_{r_{ik}^2}(u, v)  \lesssim \left( \frac{|u|^2}{r_{ik}^2} \frac{1}{V_i(|u|)} + \frac{1}{V_i(r_{ik})} \right) \, \exp\left(- c \frac{|u|^2}{r_{ik}^2}\right) \lesssim \frac{1}{V_i(|u|)},
\end{eqnarray*}
by the trivial inequality $r^{\alpha} \, e^{-r} \lesssim_{\alpha} 1$ for any $\alpha > 0$ and all $r > 0$, and the fact that $V_i(|u|)/V_i(r_{ik})$ has at most polynomial growth w.r.t. $|u|/r_{ik}$ from the doubling property. Then using the doubling property again,  \eqref{isolate-balls} gets that $h_{r_{ik}^2}(u, v) \lesssim V_i(R_i)^{-1}$.

For the opposite case $r_{ik}<1$,  similarly, using the symmetric property of the heat kernel, from \eqref{stu} and \eqref{isolate-balls},  we can write
\begin{align*}
h_{r_{ik}^2}(u, v) = h_{r_{ik}^2}(v, u)& \lesssim \frac{1}{V(v, r_{ik})} \, e^{-c \, \frac{d(v, u)^2}{r_{ik}^2}} \lesssim \frac{1}{V(v, r_{ik})} \, e^{-c' \, \frac{d(v, o_i)^2}{r_{ik}^2}} \, e^{- c' \,  \frac{R_i^2}{r_{ik}^2}} \lesssim \frac{1}{V_i(r_{ik})} \, e^{- c' \,  \frac{R_i^2}{r_{ik}^2}} \lesssim \frac{1}{V_i(R_i)},
\end{align*}
where the penultimate inequality follows from Lemma \ref{vg} (iii) and a standard trick of doubling property (cf. e.g. the proof of \eqref{BE1}).

In conclusion, we yield \eqref{nSGn}, which together with $\left\| e^{-r_{ik}^2 \L} \right\|_{L^\infty(M \setminus E_i)\to L^\infty(B_{ik})} \le 1$ implies via the  Riesz-Thorin interpolation theorem that
$$\left\| e^{-r_{ik}^2 \L} \right\|_{L^2(M \setminus E_i) \to L^\infty(B_{ik})} \lesssim \frac{1}{\sqrt{V_i(R_i)}}.$$

Hence, for any $1 \le i \le \ell$ and $\|\psi\|_2 \le 1$, we have that
\begin{eqnarray} \label{nEn}
\left|\int_{\bigcup_{0 \le j\neq i \le \ell}E_j} \left( \sum_k e^{-r_{ik}^2\L}b_{ik} \right) \, \psi\,d\mu\right|&&\le
\int_{M} \sum_k |b_{ik}| \, e^{-r_{ik}^2\L}(|\psi\chi_{M \setminus E_i}|) \,d\mu \nonumber\\
&&\lesssim
\sum_{k} \int_{B_{ik}} \frac{|b_{ik}|}{\sqrt{V_i(R_i)}}\,d\mu
\lesssim \lambda \, \frac{\sum_k\mu(B_{ik})}{\sqrt{V_i(R_i)}},
\end{eqnarray}
where we have used, in the last step, the H\"older inequality and \eqref{cz-2}.
Combining this with the fact that $\sum_k\mu(B_{ik})\lesssim \mu(F_i^{(R_i)})\sim V_i(R_i)$ (see \eqref{cz-measure}, \eqref{fi} and \eqref{brN}), we can conclude that
\begin{align*}
\cg_{22} \lesssim  \frac{1}{\lambda^2} \left(\lambda \, \frac{\sum_k\mu(B_{ik})}{\sqrt{V_i(R_i)}}\right)^2 \lesssim \sum_k\mu(B_{ik})\lesssim \frac{\|f\|_p^p}{\lambda^p},
\end{align*}
which finishes the proof.
\end{proof}

We now turn to $\cg_1$. Recall that  $p_0\in (1,2]$ is given in \eqref{set-index-1}.
\begin{prop}\label{est-away-center-offdiagonal-2}
Let $1<p<p_0$. It holds {uniformly for all $f \in L^p(M)$ and $0 < \lambda \ll \|f\|_{p}$} that
\begin{eqnarray*}
\cg_1 = \mu\left(\left\{x\in \bigcup_{1 \le j\neq i \le \ell} \left( E_j\setminus F_j^{(R_j)} \right):\left|T\left(\sum_k (1-e^{-r_{ik}^2\L}) b_{ik}\right)\right|(x) > \frac{\lambda}{2} \right\}\right)\lesssim \frac{\|f\|_p^p}{\lambda^p}.
\end{eqnarray*}
\end{prop}
\begin{proof}
By the Chebyshev inequality, we can write
\begin{eqnarray*}
\cg_1 &&\lesssim \frac{1}{\lambda^p}\int_{\bigcup_{1 \le j\neq i \le \ell} (E_j\setminus F_j^{(R_j)})} \left|T\left(
\sum_k (1-e^{-r_{ik}^2\L}) b_{ik}\right)\right|^p\,d\mu\\
&&\lesssim \frac{1}{\lambda^p}\int_{\bigcup_{1 \le j\neq i \le \ell} (E_j\setminus F_j^{(R_j)})} \left|\sum_k \left(\int_1^\infty\int_0^{r_{ik}^2}|\nabla\L e^{-(s+t)\L} b_{ik}|\frac{\,ds\,dt}{\sqrt t}\right)\right|^p\,d\mu,
\end{eqnarray*}
from the definition of $T$. Then the Minkowski's integral inequality yields that
\begin{align*}
\cg_1 &\lesssim \frac{1}{\lambda^p}\left[\sum_k\int_1^\infty\int_0^{r_{ik}^2}\left(\int_{\bigcup_{1 \le j\neq i \le \ell} E_j\setminus F_j^{(R_j)}}|\nabla \L e^{-(t+s)\L}b_{ik}|^p\,d\mu\right)^{1/p}\frac{\,ds\,dt}{\sqrt t}\right]^p\\
&\lesssim  \frac{1}{\lambda^p}\left[\sum_k\int_1^\infty\int_0^{r_{ik}^2}\|b_{ik}\|_p e^{-c \, \frac{\dist(B_{ik},E_0)^2}{t+s}} \frac{\,ds\,dt}{(t+s)^{3/2}\sqrt t}\right]^p,
\end{align*}
because of Corollary \ref{davies-gaffney-com}. Using \eqref{isolate-balls-2} and \eqref{isolate-balls-1}, the double integral can be controlled by
\begin{eqnarray*}
&&\int_1^\infty\int_0^{r_{ik}^2} e^{-c \frac{d(x_{ik},o_i)^2}{t+s}} \frac{\,ds\,dt}{(t+s)^{3/2}\sqrt t}\\
&&\quad\lesssim \int_1^{d(x_{ik},o_i)^2}\int_0^{r_{ik}^2}e^{- c \frac{d(x_{ik},o_i)^2}{t+s}} \frac{\,ds\,dt}{(t+s)^{3/2}\sqrt t}+\int_{d(x_{ik},o_i)^2}^\infty\int_0^{r_{ik}^2}e^{-c \frac{d(x_{ik},o_i)^2}{t+s}} \frac{\,ds\,dt}{(t+s)^{3/2}\sqrt t}\\
&&\quad\lesssim \int_1^{d(x_{ik},o_i)^2} \frac{dt}{\sqrt t} \int_0^{r_{ik}^2} d(x_{ik},o_i)^{-3} \, ds + \int_{d(x_{ik},o_i)^2}^\infty \frac{dt}{t^2} \int_0^{r_{ik}^2} \, ds \lesssim \frac{r_{ik}^2}{d(x_{ik},o_i)^2}.
\end{eqnarray*}

As a consequence,
\begin{eqnarray*}
\cg_1 \lesssim \frac{1}{\lambda^p}\left[\sum_k \frac{r_{ik}^2}{d(x_{ik},o_i)^2}\|b_{ik}\|_p \right]^p.
\end{eqnarray*}
Next notice that $p'>N_\infty/2$ provided $1< p < p_0$. From \eqref{isolate-balls-2} and \eqref{c2}, we have that
$$\left(\frac{r_{ik}^2}{d(x_{ik},o_i)^2}\right)^{p'}\le \left(\frac{r_{ik}}{d(x_{ik},o_i)}\right)^{N_i}\lesssim \frac{V_i(x_{ik},r_{ik})}{V_{i}(x_{ik},d(x_{ik},o_i))}\lesssim
\frac{\mu(B_{ik})}{V_{i}(o_i,d(x_{ik},o_i))}\lesssim  \frac{\mu(B_{ik})}{V_i(R_i)},$$
by means of the doubling property. Therefore, using \eqref{cz-measure}, \eqref{fi} and \eqref{brN}, we obtain that
\begin{equation*}
\sum_k \left(\frac{r_{ik}^2}{d(x_{ik},o_i)^2}\right)^{p'}\lesssim
\frac{\sum_k\mu(B_{ik})}{V_i(R_i)}\lesssim 1.
\end{equation*}
By the H\"older inequality, we conclude that
\begin{eqnarray}\label{sum-measure}
\cg_1 \lesssim \frac{1}{\lambda^p}\left[\sum_k \frac{r_{ik}^2}{d(x_{ik},o_i)^2}\|b_{ik}\|_p \right]^p&&\lesssim \frac{1}{\lambda^p} \left[\sum_k \left(\frac{r_{ik}^2}{d(x_{ik},o_i)^2}\right)^{p'}\right]^{p/p'}\left[\sum_{k}\|b_{ik}\|_p^p\right]\nonumber\\
&&\lesssim \frac{1}{\lambda^p}\left[\lambda^p\sum_k\mu(B_{ik})\right]\lesssim \frac{1}{\lambda^p}\|f\chi_{E_i\setminus F_i^{(R_i)}}\|_p^p\lesssim \frac{\|f\|_p^p}{\lambda^p},
\end{eqnarray}
where we have used the properties of Calder\'on-Zygmund decomposition \eqref{cz-2} and \eqref{cz-3}.

This completes the proof.
\end{proof}

\subsubsection{Estimate of the `bad' parts on the diagonal}
\hskip\parindent It remains to estimate
 $$\cs := \mu\left(\left\{x\in  E_i\setminus F_i^{(R_i)}:\left|T\left(\sum_k  b_{ik}\right)\right|>\lambda\right\}\right).$$
Recall that
$$T = \frac{1}{\sqrt{\pi}} \int_1^\infty \nabla e^{-t\L}\frac{\,dt}{\sqrt t}, \qquad \nabla \L^{-1/2}=\frac{1}{\sqrt{\pi}} \int_0^\infty \nabla e^{-s\L}\frac{\,ds}{\sqrt s}.$$
We then have
\begin{align*}
\cs &\le \mu\left(\left\{x\in E_i\setminus F_i^{(R_i)}:\left|\nabla \L^{-1/2}\left(\sum_k b_{ik}\right)\right|>\frac{1}{2} \lambda\right\}\right)\\
&+ \mu\left(\left\{x\in E_i\setminus F_i^{(R_i)}:\left|\int_0^1\nabla e^{-t\L}\left(\sum_k b_{ik}\right)\frac{\,dt}{\sqrt t}\right|> \frac{1}{2} \lambda \right\}\right).
\end{align*}

Let us begin with the following:
\begin{lem}\label{est-away-center-diagonal-1}
Let $1<p<2$. We have uniformly for all $f  \in L^p(M)$  and all $0 < \lambda <\|f\|_{p}$ that
\begin{eqnarray*}
 \cs_1(\lambda) =  \mu\left(\left\{x\in E_i\setminus F_i^{(R_i)}:\left|\int_0^1\nabla e^{-t\L}\left(\sum_k b_{ik}\right)\frac{\,dt}{\sqrt t}\right|> \lambda\right\}\right) \lesssim_p \frac{\|f\|_p^p}{\lambda^p}.
\end{eqnarray*}
\end{lem}
\begin{proof}
Recall that the operator $\int_0^1\nabla e^{-t\L}\frac{\,dt}{\sqrt t}$
is bounded on $L^p(M)$, cf. Lemma \ref{local-part}. The Chebyshev inequality implies that
\begin{eqnarray*}
\cs_1(\lambda) \lesssim_p \frac{1}{\lambda^p}\left\|\sum_kb_{ik}\right\|_{p}^p \lesssim \frac{1}{\lambda^p} \sum_k \left\|b_{ik}\right\|_{p}^p \lesssim \frac{\|f\|_p^p}{\lambda^p},
\end{eqnarray*}
where we have used in the second inequality the fact that the supports of $\{b_{ik}\}_{k}$ are of bounded overlap, \eqref{cz-2} and \eqref{cz-3} in the last one.
\end{proof}
For the remaining term, we write
\begin{eqnarray*}
&& \mu\left(\left\{x\in E_i\setminus F_i^{(R_i)}:\left|\nabla \L^{-1/2}\left(\sum_k b_{ik}\right)\right| > \lambda\right\}\right)\\
&&\quad\le  \mu\left(\left\{x\in E_i\setminus F_i^{(R_i)}:\left|\nabla \L^{-1/2}\left(\sum_k e^{-r_{ik}^2\L} b_{ik}\right)\right|>\frac{1}{2} \lambda\right\}\right)\\
&&\quad\quad+\mu\left(\left\{x\in E_i\setminus F_i^{(R_i)}:\left|\nabla \L^{-1/2}\left(\sum_k (1-e^{-r_{ik}^2\L}) b_{ik}\right)\right|> \frac{1}{2} \lambda\right\}\right).
\end{eqnarray*}

\begin{lem}\label{est-away-center-diagonal-2}
Let $1<p<2$.  It holds that for all $f  \in L^p(M)$  and all $0 < \lambda < \|f\|_{p}$
$$\cs_{21}(\lambda) := \mu\left(\left\{x\in E_i\setminus F_i^{(R_i)}:\left|\nabla \L^{-1/2}\left(\sum_k e^{-r_{ik}^2\L} b_{ik}\right)\right| > \lambda\right\}\right) \lesssim \frac{\|f\|_p^p}{\lambda^p}.  $$
\end{lem}
\begin{proof}
By the natural $L^2$-boundedness of $\nabla \L^{-1/2}$, we obtain that
\begin{eqnarray*}
\cs_{21}(\lambda) \lesssim \frac{1}{\lambda^2}\left\| \sum_k e^{-r_{ik}^2\L} b_{ik}\right\|_2^2  \lesssim \frac{\|f\|_p^p}{\lambda^p},
\end{eqnarray*}
by recalling that the last inequality is obtained in the proof of Proposition \ref{est-away-center-offdiagonal-1}.
\end{proof}

\begin{prop}\label{est-away-center-diagonal-3}
Let $1<p<p_0$.  We have that for all $f  \in L^p(M)$  and all $0 < \lambda \ll \|f\|_{p}$
$$\cs_{22}(\lambda) := \mu\left(\left\{x\in E_i\setminus F_i^{(R_i)}:\left|\nabla \L^{-1/2}\left(\sum_k (1-e^{-r_{ik}^2\L}) b_{ik}\right)\right|> \lambda\right\}\right) \lesssim \frac{\|f\|_p^p}{\lambda^p}. $$
\end{prop}
\begin{proof}
As in \cite[pp. 1160-1161]{cd99},  using \eqref{isolate-balls-2},  we write
\begin{eqnarray*}
\sqrt{\pi} \nabla \L^{-1/2}(1-e^{-r_{ik}^2\L}) &&=\int_{0}^\infty\left[\frac{1}{\sqrt t}-\frac{\chi_{\{t>r_{ik}^2\}}}{\sqrt{t-r_{ik}^2}}\right]\nabla e^{-t\L}\,dt\\
&&=\int_{0}^{R_i^2}\left[\frac{1}{\sqrt t}-\frac{\chi_{\{t>r_{ik}^2\}}}{\sqrt{t-r_{ik}^2}}\right]\nabla e^{-t\L}\,dt+\int_{R_i^2}^{4d(x_{ik},o_i)^2}\cdots+\int_{4d(x_{ik},o_i)^2}^\infty\cdots\\
&&=:T_{ik1}+T_{ik2}+T_{ik3}.
\end{eqnarray*}
It holds then
\begin{eqnarray}
\cs_{22}(\lambda) \le \mu(\cup_k 2B_{ik})+\sum_{j=1}^3 \cs_{22j},
\end{eqnarray}
with
\begin{align*}
\cs_{22j} := \mu\left(\left\{x\in E_i \setminus \left(F_i^{(R_i)} \cup (\cup_k 2B_{ik}) \right): \left|\sum_k T_{ikj}b_{ik}\right| > \frac{1}{3} \lambda\right\}\right), \quad 1 \le j \le 3.
\end{align*}

The estimation of $\mu(\cup_k 2B_{ik})$ is standard, via the doubling property and \eqref{cz-measure}.

Next, we estime $\cs_{221}$. Recall that \eqref{isolate-balls} says $B_{ik}\subset E_i \setminus \{z\in E_i:\,\dist(z,E_0)>R_i/2\} = E_i \setminus F_i^{(R_i/4)}$. Then for $t\le R_i^2$ and  $y\in B_{ik}$, by Proposition \ref{est-integral-diagonal}, it holds
$$\int_{E_i\setminus (F_i^{(R_i)}\cup 2B_{ik})}|\nabla h_t(x,y)|\,d\mu(x)\lesssim \frac{1}{\sqrt t} e^{- c \frac{r_{ik}^2}{t}}.$$
Therefore, following the argument in \cite[p. 1161]{cd99}, we deduce that
\begin{align*}
\cs_{221} &\lesssim \frac{1}{\lambda}\sum_k\int_{0}^{R_i^2}\int_{E_i\setminus (F_i^{(R_i)} \cup 2B_{ik})}
\left|\frac{1}{\sqrt t}-\frac{\chi_{\{t>r_{ik}^2\}}}{\sqrt{t-r_{ik}^2}}\right| |\nabla e^{-t\L}b_{ik}|\,d\mu\,dt\\
&\lesssim \frac{1}{\lambda}\sum_k\int_{0}^{R_i^2} \frac{1}{\sqrt t} e^{-c \frac{r_{ik}^2}{t}}\left|\frac{1}{\sqrt t}-\frac{\chi_{\{t>r_{ik}^2\}}}{\sqrt{t-r_{ik}^2}}\right|\|b_{ik}\|_1\,dt\\
&\lesssim \frac{1}{\lambda}\sum_k\|b_{ik}\|_1 \left(\int_{0}^{r_{ik}^2} \frac{1}{\sqrt t} e^{-c \frac{r_{ik}^2}{t}}\frac{1}{\sqrt t}\,dt+\int_{r_{ik}^2}^\infty \frac{1}{\sqrt t}\left|\frac{1}{\sqrt t}-\frac{1}{\sqrt{t-r_{ik}^2}}\right|\,dt\right) \\
&\lesssim \frac{1}{\lambda}\sum_k\|b_{ik}\|_1
\lesssim \frac{\|f\|_p^p}{\lambda^p},
\end{align*}
where we have used \eqref{nEn} in the last inequality.

Aiming now at $\cs_{223}$. From the mapping property of $\nabla e^{-t\L}$ (Proposition \ref{mapping-gradient-heat}), we get that
\begin{align*}
\cs_{223} &\lesssim  \frac{1}{\lambda^p}\int_{E_i\setminus F_i^{(R_i)}} \left|\sum_k T_{ik3}b_{ik}\right|^p\,d\mu\\
&\lesssim  \frac{1}{\lambda^p}\int_{E_i\setminus F_i^{(R_i)}} \left|\int_0^\infty\sum_k  \chi_{(4d(x_{ik},o_i)^2,\infty)}(t) \left|\frac{1}{\sqrt t}-\frac{\chi_{\{t>r_{ik}^2\}}}{\sqrt{t-r_{ik}^2}}\right| |\nabla  e^{-t\L}b_{ik}|\,dt \right|^p\,d\mu\\
&\lesssim \frac{1}{\lambda^p}\left[\sum_k\int_{4d(x_{ik},o_i)^2}^\infty \left|\frac{1}{\sqrt t}-\frac{\chi_{\{t>r_{ik}^2\}}}{\sqrt{t-r_{ik}^2}}\right| \left(\int_{E_i\setminus F_i^{(R_i)}}|\nabla e^{-t\L}b_{ik}|^p\,d\mu\right)^{1/p}\,dt\right]^p\\
&\lesssim  \frac{1}{\lambda^p}\left[\sum_k \|b_{ik}\|_p \int_{4d(x_{ik},o_i)^2}^\infty \left|\frac{1}{\sqrt t}-\frac{1}{\sqrt{t-r_{ik}^2}}\right|  \frac{\,dt}{\sqrt t}\right]^p,
\end{align*}
and in view of \eqref{isolate-balls-2}, the last integral here equals
\begin{align*}
r_{ik}^2 \, \int_{4d(x_{ik},o_i)^2}^\infty  \frac{1}{\sqrt t+\sqrt {t-r_{ik}^2}} \, \frac{1}{\sqrt{t-r_{ik}^2} \sqrt t}  \frac{\,dt}{\sqrt t}
\lesssim \frac{r_{ik}^2}{d(x_{ik},o_i)^2}.
\end{align*}
Consequently,  by \eqref{sum-measure}, we obtain that
\begin{align*}
\cs_{223} \lesssim \frac{1}{\lambda^p} \left[\sum_k \frac{r_{ik}^2}{d(x_{ik},o_i)^2}\|b_{ik}\|_p\right]^p \lesssim  \frac{\|f\|_p^p}{\lambda^p}.
\end{align*}

In conclusion, it remains to prove the lemma below.
\end{proof}

\begin{lem}
It holds that $\cs_{222} \lesssim \lambda^{-p}  \, \| f \|_p$ provided $0 < \lambda \ll \| f \|_p$.
\end{lem}
\begin{proof}
By \eqref{isolate-balls-1} and \eqref{isolate-balls-2},
\begin{equation*}
\frac{1}{2}(R_i+r_{ik})\le \dist(B_{ik},E_0)<d(x_{ik},o_i)\le 1+r_{ik}+\dist(B_{ik},E_0)\le 3\dist(B_{ik},E_0),
\end{equation*}
we split the set $E_i\setminus (F_i^{(R_i)} \cup (\cup_k 2B_{ik}))$ into
$$G_{1k}:=\{x\in E_i\setminus (F_i^{(R_i)} \cup (\cup_k 2B_{ik})): d(x,E_0)\le \frac 12\dist(B_{ik},E_0)\}$$
and
$$G_{2k}:=\{x\in E_i\setminus (F_i^{(R_i)} \cup (\cup_k 2B_{ik})): d(x,E_0)> \frac 12\dist(B_{ik},E_0)\}.$$
It holds then
\begin{gather}
\dist(G_{1k},B_{ik})\ge \frac 12\dist(B_{ik},E_0)\ge \frac 16 d(x_{ik},o_i), \label{nJn1} \\
\dist(G_{2k},E_0)\ge \frac 12\dist(B_{ik},E_0)\ge \frac 16d(x_{ik},o_i). \label{nJn2}
\end{gather}
By Proposition \ref{est-integral-diagonal}, for $t \le 4 d(x_{ik},o_i)^2$, $y\in
\{x\in E_i:\,\dist(x,E_0)\ge \frac 16d(x_{ik},o_i)\}\cap B_{ik}$, it holds
\begin{equation}\label{gradient-heat}
\int_{\{x\in E_i:\,\dist(x,E_0)\ge \frac 16d(x_{ik},o_i)\} \setminus ( 2B_{ik})}
|\nabla h_t(x,y)|\,d\mu(x)\lesssim \frac{1}{\sqrt t}e^{-\frac{r_{ik}^2}{ct}}.
\end{equation}

Notice that
\begin{align*}
\cs_{222} &=\mu\left(\left\{x\in E_i\setminus (F_i^{(R_i)} \cup (\cup_k 2B_{ik})): \left| \sum_k (\chi_{G_{1k}}+\chi_{G_{2k}})T_{ik2}b_{ik} \right| > \frac{\lambda}{3} \right\}\right) \\
&\lesssim  \frac{1}{\lambda^p} \int_{E_i\setminus F_i^{(R_i)}} \left| \sum_k \chi_{G_{1k}}T_{ik2}b_{ik} \right|^p\,d\mu +
\frac{1}{\lambda}\int_{E_i\setminus F_i^{(R_i)}}  \left|\sum_k \chi_{G_{2k}}T_{ik2}b_{ik}\right|\,d\mu =: \cf_1 + \cf_2.
\end{align*}

We begin with $\cf_1$. We argue as in
the estimation of $\cs_{223}$, by using the Davies-Gaffney estimate (cf. Corollary \ref{davies-operators}) instead of Proposition \ref{mapping-gradient-heat}, and conclude that
\begin{align*}
\cf_1 &\sim \frac{1}{\lambda^p}\int_{E_i\setminus F_i^{(R_i)}} \left|\int_0^\infty \sum_k \chi_{(R_i^2,4d(x_{ik},o_i)^2)}(t)\left|\frac{1}{\sqrt t}-\frac{\chi_{\{t>r_{ik}^2\}}}{\sqrt{t-r_{ik}^2}}\right| \, | \chi_{G_{1k}}\nabla e^{-t\L}b_{ik} | \,dt\right|^p\,d\mu\\
&\lesssim \frac{1}{\lambda^p}\left[\sum_k\int_{R_i^2}^{4d(x_{ik},o_i)^2}\left|\frac{1}{\sqrt t}-\frac{\chi_{\{t>r_{ik}^2\}}}{\sqrt{t-r_{ik}^2}}\right| \left(\int_{G_{1k}}|\nabla  e^{-t\L}b_{ik}|^p\,d\mu\right)^{1/p}\,dt\right]^p\\
&\lesssim  \frac{1}{\lambda^p} \left[ \sum_k \|b_{ik}\|_p \int_{R_i^2}^{4d(x_{ik},o_i)^2} e^{-c \frac{\dist(G_{1k},B_{ik})^2}{t}} \left|\frac{1}{\sqrt t}-\frac{\chi_{\{t>r_{ik}^2\}}}{\sqrt{t-r_{ik}^2}}\right|  \frac{\,dt}{\sqrt t} \right]^p,
\end{align*}
and from \eqref{nJn1}, the last integral here can be controlled by
\begin{align*}
\int_{R_i^2}^{4d(x_{ik},o_i)^2} e^{-c \frac{d(x_{ik},o_i)^2}{t}} \left|\frac{1}{\sqrt t}-\frac{\chi_{\{t>r_{ik}^2\}}}{\sqrt{t-r_{ik}^2}}\right|  \frac{\,dt}{\sqrt t}
&\lesssim \int_0^{r_{ik}^2} \frac{dt}{d(x_{ik}, o_i)^2} + \int_{r_{ik}^2}^{4 d(x_{ik},o_i)^2} \frac{r_{ik}^2}{d(x_{ik}, o_i)^3} \frac{dt}{\sqrt{t - r_{ik}^2}} \\
&\lesssim \frac{r_{ik}^2}{d(x_{ik}, o_i)^2}.
\end{align*}
Inserting it in the last inequality and using \eqref{sum-measure}, we get $\cf_1 \lesssim \lambda^{-p}  \, \| f \|_p$.

For $\cf_2$, by using \eqref{gradient-heat}, the argument used in the estimation of $\cs_{221}$ implies that
\begin{align*}
\cf_2 &\lesssim \frac{1}{\lambda}\sum_k\int_{R_i^2}^{4d(x_{ik},o_i)^2}\int_{G_{2k}}
\left|\frac{1}{\sqrt t}-\frac{\chi_{\{t>r_{ik}^2\}}}{\sqrt{t-r_{ik}^2}}\right| \, | \nabla e^{-t\L}b_{ik} | \,d\mu\,dt \\
&\lesssim \frac{1}{\lambda} \sum_k \| b_{ik} \|_1 \int_{R_i^2}^{4d(x_{ik},o_i)^2} \frac{1}{\sqrt t} e^{- c \frac{r_{ik}^2}{t}} \left|\frac{1}{\sqrt t}-\frac{\chi_{\{t>r_{ik}^2\}}}{\sqrt{t-r_{ik}^2}} \right| \,dt
\lesssim \frac{\|f\|_p^p}{\lambda^p}.
\end{align*}

This completes the proof.
\end{proof}

\subsection{The part around the center}\hskip\parindent
For the part around the center, we have
\begin{prop}\label{est-around-center}
Let $1<p<2$. It holds that
\begin{eqnarray*}
\ca = \sum_{i=1}^\ell
\mu\left(\left\{x \in \bigcup_{j=1}^\ell \left( E_j\setminus F_j^{(R_j)} \right): |T(f\chi_{F_i^{(R_i)}})|>\lambda\right\}\right)\lesssim \frac{\|f\|_p^p}{\lambda^p},  \qquad \forall \lambda \ll \|f\|_{p}, \  f \in L^p.
\end{eqnarray*}
\end{prop}
\begin{proof}
For each $j$, we write
\begin{eqnarray}\label{key-est-around-neck}
&&\mu\left(\left\{x\in  E_j\setminus F_j^{(R_j)}:\, |T(f\chi_{F_i^{(R_i)}})(x)|> \lambda\right\}\right)\nonumber\\
&&\quad\le
\mu\left(\left\{x\in  E_j\setminus F_j^{(R_j)}:\, \left| T\left( f\chi_{F_i^{(R_i)}} - e^{-R_j^2\L} (f \chi_{F_i^{(R_i)}}) \right)(x) \right| > \frac{1}{2} \lambda\right\}\right)\nonumber\\
&&\quad\quad+\mu\left(\left\{x\in  E_j\setminus F_j^{(R_j)}:\, \left| T\left( e^{-R_j^2\L} (f \chi_{F_i^{(R_i)}}) \right)(x) \right| > \frac{1}{2} \lambda\right\}\right) {  := \K_{j, i, 1} + \K_{j, i, 2}.}
\end{eqnarray}

{\bf Estimation of $\K_{j, i, 2}$:}
Using the $L^2$-boundedness of $T$ (Lemma \ref{local-part}), we deduce that
\begin{eqnarray} \label{K2e}
\K_{j, i, 2} &&\lesssim \frac{1}{\lambda^2}\int_{E_j\setminus F_j^{(R_j)}} \left|T \left( e^{-R_j^2\L} (f\chi_{F_i^{(R_i)}}) \right) \right|^2 \,d\mu \nonumber \\
&&\lesssim \frac{1}{\lambda^2}\left[\int_1^\infty \left\|\nabla e^{-(t+R_j^2)\L}(f\chi_{F_i^{(R_i)}})\right\|_{L^2(E_j\setminus F_j^{(R_j)})}\frac{\,dt}{\sqrt t}\right]^2.
\end{eqnarray}

Let $\psi_j$ be a Lipschitz function  such that $\psi_j\equiv 1$ on $E_j\setminus F_j^{(R_j)}$ with $\supp \psi_j\subset \{x:\,\dist(x,E_j\setminus F_j^{(R_j)})<R_j\}$ and $|\nabla \psi_j|\le C/R_j$.
Note that it holds
$$\dist(\supp \psi_j,\,E_0)>R_j.$$
We use a {Caccioppoli's inequality} as follows,
\begin{eqnarray*}
&&\int_{M}|\nabla e^{-s\L}g|^2 \, \psi_j^2\,d\mu\\
&&\quad=\int_M \nabla e^{-s\L}g \cdot \nabla  (\psi_j^2e^{-s\L}g)\,d\mu-2\int_M \left( \nabla e^{-s\L}g \cdot \nabla  \psi_j \right) \, \psi_j e^{-s\L}g\,d\mu\\
&&\quad\le \int_M (\L e^{-s\L}g) \,  (\psi_j^2e^{-s\L}g)\,d\mu+\frac{1}{2}\int_{M}|\nabla e^{-s\L}g|^2\psi_j^2\,d\mu+2\int_{\{x\in E_j:\,R_j\le |x|\le 2R_j\}}|\nabla \psi_j|^2 (e^{-s\L}g)^2\,d\mu
\end{eqnarray*}
This implies
\begin{eqnarray*}
&&\int_{E_j\setminus F_j^{(R_j)}} |\nabla e^{-s\L}g|^2\,d\mu \le \int_{M}|\nabla e^{-s\L}g|^2 \, \psi_j^2\,d\mu\\
&&\quad\le 2\left(\int_M |\L e^{-s\L}g|^2\psi_j^2\,d\mu\right)^{1/2} \, \left(\int_M|\psi_j \, e^{-s\L}g|^2\,d\mu\right)^{1/2}+\frac{C}{R_j^2}\int_{\{x\in E_j:\,R_j\le |x|\le 2R_j\}}(e^{-s\L}g)^2\,d\mu.
\end{eqnarray*}

Combining this with Propositions \ref{heat-map} and \ref{time-heat-map}, we deduce that for $1<p< 2$,
\begin{eqnarray*}
&&\left\|\nabla e^{-(t+R_j^2)\L}(f\chi_{F_i^{(R_i)}})\right\|^2_{L^2(E_j\setminus F_j^{(R_j)})}\\
&&
\quad\lesssim \frac{\|(f\chi_{F_i^{(R_i)}})\|_{p}^2}{t+R_j^2}\left(\frac{R_j}{\sqrt{ t+R_j^2}}\right) \mu(F_i^{(R_i)})^{1-\frac 2p}+\frac{\|e^{-(t+R_j^2)\L}(f\chi_{F_i^{(R_i)}})\|_{L^\infty(E_j\setminus F_j)}^2}{R_j^2}V_j(2R_j)\\
&&
\quad\lesssim \frac{\|(f\chi_{F_i^{(R_i)}})\|_{p}^2}{t+R_j^2}\left(\frac{R_j}{\sqrt{t+R_j^2}}\right) \mu(F_i^{(R_i)})^{1-\frac 2p}+\frac{\|(f\chi_{F_i^{(R_i)}})\|_{p}^2}{R_j^2}\left(\frac{R_j^{2}}{t+R_j^2}\right)^2 \mu(F_i^{(R_i)})^{1-\frac 2p}\\
&&\quad\lesssim \frac{\|(f\chi_{F_i^{(R_i)}})\|_{p}^2}{t+R_j^2}\frac{R_j}{\sqrt{t+R_j^2}} \mu(F_i^{(R_i)})^{1-\frac 2p}.
\end{eqnarray*}

Inserting this in  \eqref{K2e}, we obtain that
\begin{eqnarray}\label{key-est-around-neck-2}
\K_{j,i,2} &&\lesssim \frac{\|(f\chi_{F_i^{(R_i)}})\|_{p}^2}{\lambda^2 } \mu(F_i^{(R_i)})^{1-\frac 2p}\left[\int_1^\infty \left(\frac{R_j}{\sqrt{t+R_j^2}}\right)^{1/2} \frac{1}{\sqrt{t+R_j^2}} \frac{\,dt}{\sqrt t} \right]^2\nonumber\\
 &&\lesssim \frac{\|(f\chi_{F_i^{(R_i)}})\|_{p}^2}{\lambda^2 } \mu(F_i^{(R_i)})^{1-\frac 2p}\nonumber\\
  &&\lesssim \frac{\|f\|_{p}^p}{\lambda^p },
\end{eqnarray}
since $\mu(F_i^{(R_i)})\sim\frac{\|f\|_p^p}{\lambda^p}$ (cf. \eqref{fi}).

{\bf Estimation of $\K_{j, i, 1}$:} For $1 \le j \neq i \le \ell$, it follows from Corollary \ref{davies-gaffney-com} that
\begin{eqnarray}\label{key-est-around-neck-3}
\K_{j, i, 1} &&\lesssim \frac{1}{\lambda^p} \int_{E_j\setminus F_j^{(R_j)}} \left| T\left( f \chi_{F_i^{(R_i)}}-e^{-R_j^2\L} (f\chi_{F_i^{(R_i)}}) \right)(x) \right|^p\,d\mu\nonumber\\
&&\lesssim \frac{1}{\lambda^p} \left[ \int_1^\infty\int_0^{R_j^2} \left\| \nabla \L e^{-(s+t)\L}(f\chi_{F_i^{(R_i)}}) \right\|_{L^p(E_j\setminus F_j^{(R_j)})}\frac{\,ds\,dt}{\sqrt t}\right]^p\nonumber\\
&&\quad\lesssim  \frac{1}{\lambda^p} \left[ \int_{R_j^2}^\infty \int_0^{R_j^2} \left(\frac{\|f\chi_{F_i^{(R_i)}}\|_p}{(s+t)^{3/2}}\right) \frac{\,ds\,dt}{\sqrt t}+\int_1^{R_j^2}\int_0^{R_j^2}\left(\frac{\|f\chi_{F_i^{(R_i)}}\|_p}{(s+t)^{3/2}}\right) e^{-c \frac{R_j^2}{s+t}} \frac{\,ds\,dt}{\sqrt t} \right]^p\nonumber\\
&&\lesssim \frac{\|f\chi_{F_i^{(R_i)}}\|_p^p}{\lambda^p}.
\end{eqnarray}

For the opposite case $j=i$, it suffices to slightly modify the above proof. More precisely, set
$$\widetilde F_i:=\{x\in E_i:\,\dist(x,E_0)\le 4R_i\}.$$
It holds then
$$\mu(\widetilde F_i)\sim V_i(R_i)\sim \lambda^{-p}\|f\|_{p}^p, \qquad \dist(E_i\setminus \widetilde F_i,\,F_i^{(R_i)})\sim R_i. $$

Now we write
$$\K_{i, i, 1} \le \mu(\widetilde F_i) + \mu\left(\left\{x\in  E_i\setminus \widetilde F_i:\,  \left| T\left( f \chi_{F_i^{(R_i)}} - e^{-R_i^2\L} (f \chi_{F_i^{(R_i)}}) \right)(x) \right| > 2^{-1} \lambda\right\}\right).$$
Repeat the argument in the case $j \neq i$, one obtains immediately $\K_{i, i, 1} \lesssim \frac{\|f\|_p^p}{\lambda^p}$.

This together with \eqref{key-est-around-neck-2} and \eqref{key-est-around-neck-3} gives the desired estimate and completes the proof.
\end{proof}

\subsection{Completion of the proof}
\hskip\parindent
Combining the estimates from the previous three subsection, i.e., the estimates on the center (Lemma \ref{est-center}), the estimates
on the part away from the center (Proposition \ref{est-away-center}), and the estimates near the center (Proposition \ref{est-around-center}), we finally conclude that for $\lambda<\|f\|_p$,
\begin{eqnarray*}
&&\mu\left(\left\{x\in \bigcup_{j=1}^\ell \left( E_j\setminus F_j^{(R_j)} \right):|Tf|>(2\ell+1)\lambda\right\}\right)\\
&&\quad\le \mu\left(\left\{x\in \bigcup_{j=1}^\ell \left( E_j\setminus F_j^{(R_j)} \right):|T(f\chi_{E_0})|>\lambda\right\}\right)+\sum_{i=1}^\ell
\mu\left(\left\{x\in \bigcup_{j=1}^\ell \left( E_j\setminus F_j^{(R_j)} \right):|T(f\chi_{F_i^{(R_i)}})|>\lambda\right\}\right)\\
&&\quad\quad+\sum_{i=1}^\ell
\mu\left(\left\{x\in \bigcup_{j=1}^\ell \left( E_j\setminus F_j^{(R_j)} \right):|T(f\chi_{E_i\setminus F_i^{(R_i)}})|>\lambda\right\}\right)\\
&&\quad\lesssim \frac{\|f\|_p^p}{\lambda^p},
\end{eqnarray*}
for each $1<p<p_0$, where $p_0$ is as in \eqref{set-index-1}.
This together with \eqref{est-small-lambda-1} yields
\begin{eqnarray*}
&&\mu(\{x\in M:|Tf|>(2\ell+1)\lambda\})\le C \, \frac{\|f\|_p^p}{\lambda^p}.
\end{eqnarray*}
By this and \eqref{est-large-lambda}, we see that the operator $T$ is weakly $L^p$ bounded for each
$1<p<p_0$.

Recall that
$$T= \frac{1}{\sqrt{\pi}} \int_1^\infty\nabla e^{-t\L}\frac{\,dt}{\sqrt t}.$$
On the other hand, by Lemma \ref{local-part},
$$\frac{1}{\sqrt{\pi}} \int_0^1\nabla e^{-t\L}\frac{\,dt}{\sqrt t}$$
is bounded on $L^p(M)$ for all $1<p<2$.
We finally conclude that $\nabla \L^{-1/2}$ is weakly $L^p$ bounded for each $1<p<p_0$,
and hence the Riesz transform $\nabla \L^{-1/2}$
is bounded on $L^p(M)$ for all $1<p<2$ by the Marcinkiewicz interpolation theorem.
This completes the proof of Theorem \ref{main-result}.

\subsection*{Acknowledgments}
\addcontentsline{toc}{section}{Acknowledgments} \hskip\parindent
R. Jiang was partially supported by NNSF of China 11922114 \& 11671039) and NSF of Tianjin (20JCYBJC01410),
H.-Q. Li was  partially supported by NNSF of China (12271102 \& 11625102).

\noindent Ren-Jin Jiang \\
\noindent  Center for Applied Mathematics\\
\noindent Tianjin University\\
\noindent  Tianjin 300072\\
\noindent China\\
\noindent {rejiang@tju.edu.cn}

\

\noindent Hong-Quan Li\\
\noindent School of Mathematical Sciences\\
\noindent Fudan University\\
\noindent 220 Handan Road\\
\noindent Shanghai 200433 \\
\noindent People's Republic of China\\
\noindent hongquan\_li@fudan.edu.cn

\

\noindent Hai-Bo Lin\\
\noindent  College of Science\\
\noindent  China Agricultural University\\
\noindent  Beijing 100083\\
\noindent  People's Republic of China\\
\noindent  haibolincau@126.com

\end{document}